\documentstyle[12pt]{article}
\newfont{\sdbl}{msbm9}
\newfont{\dbl}{msbm10 at 12pt}
\newtheorem{th}{Theorem}
\newtheorem{defi}{Definition }
\newtheorem{prop}{Proposition}
\newtheorem{lemma}{Lemma}
\newtheorem{corol}{Corollary}
\newcommand{\eq}{=}
\newcommand{\dz}{{\mbox{\dbl Z}}}

\newcommand{\dn}{{\mbox{\dbl N}}}
\newcommand{\df}{{\mbox{\dbl F}}}
\newcommand{\sdf}{{\mbox{\sdbl F}}}
\begin{document}

\bigskip

\centerline{\Large Wild division algebras over Laurent series fields }

\bigskip

\centerline{{\large  A.B. Zheglov}
\footnote{Supported by the Graduiertenkolleg "Geometrie und Nichtlineare
Analisys" of the DFG }
\footnote{e-mail address: azheglov@mathematik.hu-berlin.de}}

\bigskip

\centerline{\bf Abstract}

\begin{quote} 
\small{In this paper we study some special classes of division algebras over
a Laurent series field  with arbitrary
residue field. We call the algebras from these classes as splittable and good splittable division algebras. It is shown that theses classes contain the group of tame division algebras. 
For the class of good division algebras 
a decomposition theorem is given. This theorem is a
generalization of the decomposition theorems for tame division algebras given
by Jacob and Wadsworth in \cite{JW}. 
For both clases we introduce a notion of a $\delta$-map and develop a technique of $\delta$-maps for division algebras from these classes. Using this technique we reprove several old well known results of Saltman and get 
the positive answer on the period-index conjecture of M.Artin: the exponent of $A$ is equal to its index for any division algebra $A$ over a $C_2$-field $F$, when $F\eq F_1((t_2))$, where $F_1$ is a $C_1$-field (see \cite{PY}, 3.4.5.).
The paper includes also some other results about splittable division algebras, which, we hope, will be useful for the further investigation of wild division algebras. }  
\end{quote}

\section{Introduction}

In this paper we study some class of division algebras over
a Laurent series field  with arbitrary residue field. Namely, we study division algebras which satisfy the 
following condition:
there exists a section $\bar{D} \hookrightarrow D$ of the residue homomorphism 
$D\rightarrow \bar{D}$,
where $D$ is a central division algebra over a complete discrete valued field $F=k((t))$. We say that these division algebras are splittable. 
If $char k=0$, all such division algebras are tame and therefore belong to the group of tame division algebras, which was carefully studied in the papers \cite{JW} and \cite{PY} even in a much more general situation of a henselian field $F$ of arbitrary characteristic. So, we consider mostly wild division algebras.

 An extensive analysis of the wild
division algebras of degree $p$ over a field $F$ with complete discrete rank 1
valuation with $char (\bar{F})\eq p$ was given by Saltman in \cite{Sa} (Tignol
in \cite{Ti} analyzed more general case of the defectless division algebras
of degree $p$ over a fild $F$ with Henselian valuation). Here we study  splittable division algebras of arbitrary index. 
This class (which is not a subgroup in $Br (F)$) contains a class of good splittable division algebras (see the definition in section 2), which posess several beautiful properties. In particular, we prove a decomposition theorem for such algebras.  This theorem is a
generalization of the decomposition theorems for tame division algebras given
by Jacob and Wadsworth in \cite{JW}. 

For arbitrary splittable division algebras we give only several assorted results, and the study of this class is far from to be complete. Nevertheless, we investigate here technical tools, which are important for the study of such algebras, and prove  a relation between the level and a higher order level for some splittable division algebras (see section 6). We hope this technique will be applied to the study of the cyclisity question for certain division algebras od degree $p^k$. 

As an application we get several results, which are partly  well known (see proposition \ref{cyclisity}) and party not. In particular, we get the positive answer on the following
conjecture: the exponent of $A$ is equal to its index for any division algebra
$A$ over a $C_2$-field $F\eq F_1((t_2))$,  where $F_1$ is a $C_1$-field.

Here is a brief overview of this paper. 

In section 2 we give a definition of splittable and good splittable division algebras and prove that all tame division algebras over $F=k((t))$ are good splittable. 

Section 3 contains the most important technical tools for the study of splittable division algebras. We define a notion of $\delta$-maps and investigate a theory of $\delta$-maps for such algebras. In this section we define also the notion of a local height, which is a possible generalization of Saltman's level. 

In section 4 we prove the period-index conjecture metioned above. This section contains also a small history of the question known to the author. We note that the proof does not use all the results from section 3. 

In section  5 we study good splittable division algebras and prove the decomposition theorem. 

In section 6 we reprove some results of Saltman about semiramified division algebras of index $p$ over $F$ using the technique from section 3. Then we define a notion of a higher order  level and prove several general properties of splittable division algebras satisfying the following condition: $Z(\bar{D})/\bar{F}$ is a simple extension. At the end of section we put several open questions.

We use the notation of \cite{JW}. We always denote by $D$ a division algebra
finite dimensional over its center $F=k((t))\eq Z(D)$. Recall that any Henselian
valuation on $F$ has a unique extension to a valuation on $D$. We denote the valuation on $F$ by $v$ and its unique extension on $D$ by $w$. 

Given a valuation $w$ on $D$, we denote by $\Gamma_{D}$ its value group, by
$V_D$ its valuation ring, by $M_D$ its maximal ideal and by $\bar{D}\eq
V_D/M_D$ its residue division ring.

By \cite{S}, p.21 one has the fundamental inequality
$$
[D:F]\ge |\Gamma_D:\Gamma_E| \cdot [\bar{D}:\bar{F}].
$$
$D$ is called defectless over $F$ if equality holds and defective otherwise.
It is known that $D$ is defectless if it has a discrete valuation of rank 1.

Jacob and Wadsworth in \cite{JW} introduced the basic homomorphism
$$
\theta_D:\Gamma_D/\Gamma_F\rightarrow Gal(Z(\bar{D})/\bar{F})
$$
induced by conjugation by elements of $D$. They showed that $\theta_D$ is
surjective and $Z(\bar{D})$ is the compositum of an abelian Galois and a
purely inseparable extension of $\bar{F}$.

We say $D$ is tame division algebra if $char (\bar{F})\eq 0$ or $char
(\bar{F})\eq q\ne 0$, $D$ is defectless over $F$, $Z(\bar{D})$ is separable
over $\bar{F}$, and $q{\not |}|ker(\theta_D)|$. We say $D$ is wild division
algebra if it is non tame.

We call a division algebra $D$ {\it inertially split} if $Z(\bar{D})$ is
separable over $\bar{F}$, the map $\theta_D$ is an isomorphism, and $D$ is
defectless over $F$.

\bigskip

{\bf Acknowledgements}

I am  grateful to Professor A. N. Parshin, Professor E.-W. Zink,   and M.
Grabitz for useful discussions and attention to my work. I am very grateful to Professor A.Wadsworth for 
carefully reading  my paper and for showing me a mistake in the very first version of this paper and to Professor V.I.Yanchevskii for valuable discussions during his visit in Berlin. Finally, I thank my wife Olga for her support and encouragement.

\section{Cohen's theorem}

Recall one definition from \cite{Zh}.
\begin{defi}
A division algebra $D$ is said to be splittable if there is a homomorphism
$\bar{D} \hookrightarrow {\cal O}_D\subset D$
that is a section of the map ${\cal  O}_D \rightarrow \bar{D}$.
\end{defi}

There is a natural question if there exists a generalization of Cohen's 
theorem, i.e. is any
central division algebra splittable or not. It is not true if a division 
algebra is not finite dimensional
over its centre, as Dubrovin's example in \cite{Zh} shows. It is not true also 
for some finite dimensional
division algebras, as the example to theorem 2.7. in \cite{Sa} shows.   
But it is true 
for tame division
algebras over complete discrete valued fields. This easily follows
from results of Jacob and Wadsworth \cite{JW} (compare with \cite{Zh}, Th.1).

\begin{th}
\label{Cohen}
Let $(F,v)$ be a valued field which is complete with respect to a discrete
rank 1 valuation $v$. Suppose $char F\eq char \bar{F}$. Let $D$ be a tame 
division
algebra with $Z(D)\eq F$ and $[D:F]< \infty$.

Then there exists a section $\bar{D}\hookrightarrow D$ of the residue
homomorphism $D\rightarrow \bar{D}$.
\end{th}

{\bf Proof.} Since $F$ is a complete field, $F$ is a Henselian field and $v$
extends uniquely to a valuation $w$ on $D$. Since $D$ is tame,
  $Z(\bar{D})/\overline{Z(D)}$ is a cyclic Galois extension.
There exists an
inertial lift $Z$ of $Z(\bar{D})$ over $F$, $Z$ is Galois over $F$,
and by classical Cohen's theorem there exists a section
$\tilde{Z}(\bar{D})\hookrightarrow Z$.

Consider the centraliser $C\eq C_D(Z)$ of $Z$ in $D$. Then we have $\bar{C}\eq
\bar{D}$.

Indeed, by Double Centraliser Theorem we have $[D:F]\eq [C:F][Z:F]$ and
$[Z:F]\eq |Gal (Z(\bar{D})/\bar{F})|$. By \cite{JW}, prop.1.7 a
homomorphism $\theta_D: \Gamma_D/\Gamma_F\rightarrow Gal(Z(\bar{D})/\bar{F})$
is surjective, so for any parameter $z$ we have
$\theta_D(w(z))\eq \sigma$, where $<\sigma >
\eq Gal (Z(\bar{D})/\bar{F})$. It is clear that $z\notin
C$. Now let $u_1,\ldots ,u_{[C:F]}$ be a $F$-basis of $C$. It is easy to see
that the elements $u_j, zu_j, \ldots , z^{n-1}u_j
$, $j\eq 1,\ldots ,[C:F]$, where $n\eq ord (\sigma )$,
the order of $\sigma $, are linearly independent, so form a basis for $D$ over
$F$. Since
$$
w(F\langle zu_j, \ldots , z^{n-1}u_j,
 j\eq 1,\ldots ,[C:F]\rangle )\cap \Gamma_C\eq 0,
$$
where $F\langle zu_j, \ldots , z^{n-1}u_j,
 j\eq 1,\ldots ,[C:F]\rangle$ denote a vector space in $D$ over $F$ generated 
by elements $u_jz^i$,
this implies that for any element $x\in D$ with
$w(x)\eq 0$ we can find elements $r_1, \ldots r_{[C:F]}\in F$ such that $x\eq
r_1u_1+\ldots +r_{[C:F]}u_{[C:F]}\mbox{\quad mod \quad} M_D$. Hence
$\bar{C}\eq \bar{D}$.

Note that $C$ is an unramified division algebra. Indeed, by \cite{JW}, th.2.8, 
th.2.9 $C$ contains a copy of the inertial lift of a maximal separable 
subfield in $\bar C$, say $\tilde C$. Then the centralizer $C_C(\tilde{C})$ 
must be a totally ramified division algebra, i.e. it is trivial and $\tilde C$ 
is a maximal subfield. So, $C$ must be unramified.

Fix an embedding $i: \bar F \hookrightarrow F$. It can be extended to the 
embedding $i':\bar Z \hookrightarrow Z$, $i'|_{\bar F}=i$ by Hensel lemma. Now 
consider the algebra $A\eq \bar{C}\otimes_{\bar Z}Z(C)$. It is easy to see 
that $A$ is an
unramified division algebra with $\bar{A}\eq \bar{C}\eq \bar{D}$.
Therefore by \cite{Az}, Th.31, $A\cong C$; so there exists a section
$\bar{D}\hookrightarrow C$.

The theorem is proved.\\
$\Box$\\

Later we will see that much more can be said about good  splittable algebras:

\begin{defi}
\label{goodsplit}
A division algebra $D$ is called good splittable if there exists a section 
$s:\bar D\hookrightarrow
D$ compatible with an embedding $i:\overline{Z(D)}\hookrightarrow Z(D)$, i.e.
$s(\overline{Z(D)})=i(\overline{Z(D)})\subset Z(D)$.

\end{defi}

It's easy to see that all tame division algebras are good splittable, because 
by Hensel lemma any embedding $\overline{Z(D)}\hookrightarrow Z(D)$ can be 
uniquely extended to any separable extension of $Z(D)$.

It is interesting to know what kind of splittable division algebras are good 
splittable. 

By theorem 3.9. in \cite{Sa} even a splittable division algebra $D$ of degree $p=char D$ is not a good splittable algebra if the level of $D$ (the notion of level we will recall in section 3, see remark to lemma \ref{svva}) is divisible by $p$. 
Nevertheless, it is an open question whether it is  true or not, for example, for division algebras with $\bar{D}=Z(\bar{D})$ such that  $\bar{D}/\bar{F}$ is a simple extension and the local height (see the definition in the same remark) is not divisible by $p$.  
We will discuss this question in section 6.

\section{Delta-maps of splittable algebras}

In this section we develop some ideas from \cite{Zh}, where some properties of 
$\delta$-maps for special kind of local skew fields were studied. Technical 
properties of $\delta$-maps play the main role in all our results. Here we 
will give a list of these properties.

Let $D$ be a finite dimensional
division algebra over a complete valued field $F\eq k((t))$. Let $w$ be a
unique extension of the valuation $v$ to $D$. We will denote by $z$ any
parameter of $D$, i.e. any element with $\langle w(z)\rangle \eq \Gamma_D$.
 Consider the ring $\dz \langle\alpha, \delta \rangle$ of noncommutative 
polinomials in
two variables. Define the map
$$
\sigma :\dz \langle\alpha , \sigma \rangle\rightarrow \dz \langle\alpha 
,\delta , \delta_i; i\ge
1\rangle ,
$$
$$
\sigma (\alpha^{a_1}\delta^{b_1}\ldots \alpha^{a_n}\delta^{b_n})\eq
\alpha^{a_1}\delta_{b_1}\ldots \delta_{b_{n-1}}\alpha^{a_n-1}\delta^{b_n},
$$
where $a_1, b_n\ge 0$, $a_i, b_j\ge 1$, $i>1$, $j<n$ for every word in $\dz
\langle\alpha , \delta \rangle$.

Let $S_i^k\in \dz \langle\alpha ,\delta \rangle$, $i\ge k$, $i\ge 1$ be 
polynomials given
by the following formula:
$$
S_i^k\eq \sum_{\tau\in S_i/G}\tau (\underbrace{\alpha\ldots
\alpha}_{i-k}\underbrace{\delta\ldots \delta}_{k}),
$$
where $S_i$ is a permutation group and $G$ is an isotropy subgroup.

\begin{lemma}{(\cite{Zh}, lemma 2)}
The polynomials $S_i^k$ satisfy the following property:
$$
S_i^i\eq \delta^i, \mbox{\quad} S_i^0\eq \alpha^i, \mbox{\quad}
S_{i+1}^{k+1}\eq \alpha S_i^{k+1}+\delta S_i^k  $$
\end{lemma}

For any splittable division algebra can be defined a notion of $\delta$-maps:

\begin{prop}{(\cite{Zh}, prop. 1,2)}
\label{ooo}
Let $D$ be a splittable division algebra.
Fix some parameter $z$ and some embedding $u: \bar{D} \hookrightarrow D$.
Then
$D$ is isomorphic to a division algebra $\bar{D}((z))$,
which is defined to be the vector space of series with multiplication defined 
by the formula
$$
zaz^{-1}\eq \alpha (a)+\delta_1(a)z+\delta_2(a)z^2+\ldots ,\mbox{\quad}
a\in \bar{D},
$$
where $\alpha :\bar{D}\rightarrow \bar{D}$ is an
automorphism and $\delta_i:\bar{D}\rightarrow \bar{D}$ are linear maps such
that the map $\delta_i$ satisfy the identity
$$\delta_i(ab)\eq \sum_{k\eq 0}^i\sigma(\delta^{i-k}\alpha)
(a)\sigma(S_i^k\alpha)
(b),\mbox{\quad} a,b\in \bar D
$$
\end{prop}

{\bf Remark} Note that the values $\sigma (S_i^k\alpha )$ and $\sigma
(\delta^{i-k}\alpha )$ belong to the subring $\dz \langle\alpha , \delta_i, 
i\ge
1\rangle$, so the formula is well defined.

Note that $\delta$-maps depend on the choice of a parameter and an embedding. 
The automorphism $\alpha$, as it easy to see, depend only on the choice of a 
parameter.
In the proposition we identify $\bar{D}$ with $u(\bar{D})$.

\begin{corol}{(\cite{Zh}, corol. 1)}
\label{formuly}
Suppose $\alpha\eq Id$. Then
$$
\delta_i(ab)\eq
\delta_i(a)b+
\sum_{k\eq 1}^{i}\delta_{i-k}(a)\sum_{(j_1,\ldots ,j_l)}
C_{i-k+1}^l\delta_{j_1}\ldots \delta_{j_l}(b),
$$
where $\delta_0\eq\alpha$ and the second sum is taken over all the vectors
$(j_1,\ldots ,j_l)$ such that $0< l\le min\{i-k+1, k\}$, $j_m\ge 1$, $\sum
j_m\eq k$.
\end{corol}

Further we will need even more general definition.

\begin{defi}
\label{maps}
In the situation of proposition \ref{ooo}
let us define maps ${}_m^{(z,u)}\delta_i: \bar{D}\rightarrow \bar{D}$, 
$m\in\dz$,
$i\in \dn$ as follows.
$$
z^maz^{-m}\eq u({}^{(z)}\alpha^m(\bar{a}))+u({}_m^{(z,u)}\delta_1(\bar{a}))z+
u({}_m^{(z,u)}\delta_2(\bar{a}))z^2+\ldots ,\mbox{\quad}
a\in u(\bar{D}).
$$

If $m\eq 0$, put ${}_m^{(z,u)}\delta_i\eq 0$.
\end{defi}

Note that ${}^{(z)}\alpha |_{Z(\bar{D})}$ does not depend on the choice of $z$.

Note that if ${}^{(z)}\alpha\eq id$, then ${}_m^{(z,u)}\delta_i\eq 0$ for 
$m\eq p^k$, where $k$
is sufficiently large, $k$  depends on $i$. Moreover, ${}_m^{(z,u)}\delta_i\eq
{}_{m+p^k}^{(z,u)}\delta_i$ for $k$ sufficiently large. We will use also the 
following
notation:
$$
{}_m^{(z,u)}\tilde{\delta_i}\eq {}_{-m}^{(z,u)}\delta_i, \mbox{\quad} 
{}_{1}^{(z,u)}\delta_i=  {}^{(z,u)}\delta_i
$$
Sometimes, we will write ${}_m\delta_i$ instead of
${}_m^{(z,u)}\delta_i$ and ${}_m^{(z,u)}\delta_i(a)$ instead of 
$u({}_m^{(z,u)}\delta_i(\bar{a}))$ whenever the context is clear.

Immediately from the definition follows

\begin{lemma}
\label{triviall}
In the situation of definition \ref{maps}
we have

(i) for $|m|>1$
$$
{}_m^{(z,u)}\delta_i(a)={}^{(z)}\alpha^{sign(m)} 
({}_{sign(m)(|m|-1)}^{(z,u)}\delta_i(a))+
{}_{sign(m)}^{(z,u)}\delta_i({}^{(z)}\alpha^{sign(m)(|m|-1)}(a))+
$$
$$
\sum_{j=1}^{i-1}{}_{sign(m)}^{(z,u)}\delta_j({}_{sign(m)(|m|-1)}^{(z,u)}\delta_{i-j}(a)),
$$
where $sign(m)=m/|m|$, $a\in \bar{D}$;

(ii) for any $m\ne 0$
$$
{}^{(z)}\alpha^{-m} ({}_{m}^{(z,u)}\delta_i)+
{}_{-m}^{(z,u)}\delta_i({}^{(z)}\alpha^{m})+
\sum_{j=1}^{i-1}{}_{-m}^{(z,u)}\delta_j({}_{m}^{(z,u)}\delta_{i-j})=0
$$
\end{lemma}

\begin{prop}
\label{flyii}

For fixed  $z,u$ from proposition \ref{ooo}  we have

(i) The maps ${}_m^{(z,u)}\delta_i$ satisfy the following identities:
$$
{}_m{\delta_i}(ab)\eq
{}_m{\delta_i}(a)\alpha^{i+m}(b)+\alpha^m(a){}_m{\delta_i}(b)+
\sum_{k\eq 1}^{i-1}{}_m{\delta_{i-k}}(a)
{}_{i-k+m}{\delta_{k}}(b)
$$

(ii) Suppose $\alpha \eq id$. Then
the maps ${}_m^{(z,u)}\delta_i$ satisfy the following identities:
$$
{}_m{\delta_i}(ab)\eq
{}_m{\delta_i}(a)b+a{}_m{\delta_i}(b)+
\sum_{k\eq 1}^{i-1}{}_m{\delta_{i-k}}(a)\sum_{(j_1,\ldots ,j_l)}
C_{i-k+m}^l\delta_{j_1}\ldots
{\delta_{j_l}}(b)
$$
where the second sum is taken over all the vectors
$(j_1,\ldots ,j_l)$ such that $0< l\le min\{i-k+m, k\}$, $j_m\ge 1$, $\sum
j_m\eq k$; $C_j^k\eq 0$ if $j\eq 0$, and $C_j^k\eq C_{j+p^q}^k$ for $q>>0$ if
$j\le 0$.
\end{prop}

{\bf Proof.} For any $a,b\in\bar{D}$ we have
$$
\alpha^m(ab)z^{m}+{}_m{\delta_1}(ab)z^{m+1}+{}_m{\delta_2}(ab)z^{m+2}+\ldots\eq
z^{m}(ab)\eq
$$
\begin{equation}
\label{(*)}
(\alpha^m(a)z^{m}+{}_m{\delta_1}(a)z^{m+1}+{}_m{\delta_2}(a)z^{m+2}+\ldots )b
\end{equation}
 If we represent the right-hand side of (\ref{(*)}) as a series with coeffitients
shifted to the left and then compare the corresponding coeffitients on the
left-hand side and right-hand side, we get some formulas for
${}_m\delta_i(ab)$. We have to prove that these formulas are the same as in
our proposition.

Let
$$z^{i+m-k}b\eq \alpha^{i+m-k}(b)z^{i+m-k}+\ldots +x'_kz^{i+m}+\ldots $$
and
$$(\alpha^m(a)z^{m}+{}_m{\delta_1}(a)z^{m+1}+{}_m{\delta_2}(a)z^{m+2}+\ldots )b
\eq \alpha^m(ab)z^m+y_{m+1}z^{m+1}+y_{m+2}z^{m+2}+\ldots
$$
Then we have
$$
y_{i+m}\eq \alpha^m(a)x'_i+\sum_{k\eq 0}^{i-1}{}_m\delta_{i-k}(a)x'_k
$$
In the proof of \cite{Zh}, prop.2 we have shown that
$$
z^{i+1-k}b\eq \alpha^{i+1-k}(b)z^{i+1-k}+\ldots +\sigma (S_i^k\alpha
)(b)z^{i+1}+\ldots
$$
 Hence $x'_k\eq \sigma (S_{i+m-1}^k\alpha )(b)$ for $k< i$. It is easy to
see that $x'_i\eq {}_m\delta_i(b)$, $x'_0\eq \alpha^{i+m}(b)$ and $\sigma
(S_{i+m-1}^k\alpha ) \eq {}_{i+m-k}\delta_k$, which proves
(i).

For $\alpha\eq id$, by corollary \ref{formuly},
$$
\sigma
(S_{i+m-1}^k\alpha )(b)\eq \sum_{(j_1,\ldots ,j_l)}
C_{i-k+m}^l\delta_{j_1}\ldots
{\delta_{j_l}}(b),
$$
where $l, j_1, \ldots ,j_l$ were defined in our proposition. This proves
(ii).\\
The proposition is proved.\\
$\Box$

\begin{lemma}{(\cite{Zh}, lemma 3 )}
\label{ozamene}

In the situation of proposition \ref{ooo}
suppose ${}^{(z,u)}_i\delta_j$ is the first map such that 
${}^{(z,u)}_i\delta_j(a)\ne 0$ for given $a\in \bar{D}$, $i\in \dz \backslash 
\{0\}$, i.e. 
${}^{(z,u)}_i\delta_1(a)\eq\ldots\eq {}^{(z,u)}_i\delta_{j-1}(a)\eq 0$, 
${}^{(z,u)}_i\delta_j(a)\ne 0$ (so we have a map $i\mapsto j(i)$). Then

(i) for $z'\eq z+u(b)z^{q+1}$, $b\in \bar{D}$ we have
${}^{(z')}\alpha^i (a)={}^{(z)}\alpha^i (a)$,
 ${}^{(z',u)}_i\delta_{k}(a)={}^{(z,u)}_i\delta_{k}(a)$ for $k<q$ and
$$
{}^{(z',u)}_i\delta_{q}(a)\eq {}^{(z,u)}_i\delta_{q}(a) 
+b'{}^{(z)}\alpha^{q+i}(a)-
{}^{(z)}\alpha^i (a)b',
$$
where $b'=\sum_{k=0}^{i-1}{}^{(z)}\alpha^k(b)$.

(ii) Suppose ${}^{(z)}\alpha^n|_{Z(\bar{D})}\eq id$, $n\ge 1$, $a\in 
Z(\bar{D})$ and \\
${}^{(z,u)}_1\delta_1({}^{(z)}\alpha^k(a))\eq\ldots\eq 
{}^{(z,u)}_1\delta_{j-1}({}^{(z)}\alpha^k(a))\eq 0$ for any $k$.

Then for $z'\eq z+u(b)z^{q+1}$, $b\in \bar{D}$ we have
${}^{(z')}\alpha^i (a)={}^{(z)}\alpha^i (a)$, ${}^{(z',u)}_i\delta_{k}(a)=
{}^{(z,u)}_i\delta_{k}(a)$ for $k<q+j$ and
$$
{}^{(z',u)}_i\delta_{q+j}(a)={}^{(z,u)}_i\delta_{q+j}(a)+
b'{}^{(z)}\alpha^q({}^{(z,u)}_i\delta_{j}(a))-{}^{(z,u)}_i\delta_{j}(a){}^{(z)}\alpha^j(b')+
$$
$$
b'\sum_{k\eq 
1}^{q}{}^{(z)}\alpha^{q-k}({}^{(z,u)}\delta_{j}({}^{(z)}\alpha^{k+i-1}(a))) -
{}^{(z,u)}_i\delta_{j}(a)\sum_{k\eq 0}^{j-1}{}^{(z)}\alpha^{k}(b),
$$
where $b'=\sum_{k=0}^{i-1}{}^{(z)}\alpha^k(b)$, if $n|q$ or ${}^{(z)}\alpha 
(a)=a$. \\

In particular, if ${}^{(z)}\alpha\eq id$ and $(i,p)=1$, then
$$
{}^{(z',u)}_i\delta_{q+j}(a)={}^{(z,u)}_i\delta_{q+j}(a)+ 
(q-j){}^{(z,u)}_i\delta_{j}(a)b
$$

(iii) for $z'\eq u(b)z$, $b\in Z(\bar{D})$, $b\ne 0$ we have
${}^{(z')}\alpha (a)={}^{(z)}\alpha (a)$, 
${}^{(z',u)}\delta_{k}(a)={}^{(z,u)}\delta_{k}(a)$ for $k<j$ and
$$
{}^{(z',u)}\delta_{j}(a)={}^{(z,u)}\delta_{j}(a){}^{(z)}\alpha (b^{-1})\cdots 
{}^{(z)}\alpha^j(b^{-1})
$$
if $i=1$.
\end{lemma}

{\bf Proof.}
 (i)
We have
$$
{z'}^ia{z'}^{-i}\eq (1+b'z^q+\ldots )z^iaz^{-i}(1+b'z^q+\ldots )^{-1}\eq 
(z^iaz^{-i}+b'z^qz^iaz^{-i}+\ldots )(1-b'z^q+\ldots )\eq
$$
$$
(z^iaz^{-i}-z^iaz^{-i}b'z^q+\ldots +b'z^qz^iaz^{-i}-\ldots )\eq
$$
$$
(z^iaz^{-i}-[{}^{(z)}\alpha^i (a)+
{}^{(z,u)}_i\delta_j(a)z^j+\ldots ]b'z^q+b'z^q[{}^{(z)}\alpha^i 
(a)+{}^{(z,u)}_i\delta_j(a)z^j+\ldots ]+\ldots
)\eq
$$
$$
(z^iaz^{-i}-[{}^{(z)}\alpha^i 
(a)b'+{}^{(z,u)}_i\delta_j(a){}^{(z)}\alpha^j(b')z^j+\ldots ]z^q+
b'{}^{(z)}\alpha^{q+i}(a)z^q+
\ldots )\eq
$$
$$
(z^iaz^{-i}+(-{}^{(z)}\alpha^i (a)b'+b'{}^{(z)}\alpha^{q+i}(a))z^q+
\ldots )\eq
{}^{(z)}\alpha^i (a)+\ldots  + {}^{(z,u)}_i\delta_{q-1}(a)z'^{q-1}+
$$
$$
({}^{(z,u)}_i\delta_{q}(a) +
b'{}^{(z)}\alpha^{q+i}(a)-{}^{(z)}\alpha^i (a)b')z'^q+ \ldots
$$

(ii)
Put $c=z'^iz^{-i}-1-b'z^{q+i}$. So, $w(c)>q+i$. Note that 
$c{}^{(z)}\alpha^k(a)={}^{(z)}\alpha^k(a)c$, since $n|q$ or ${}^{(z)}\alpha 
(a)=a$ and $a\in Z(\bar{D})$. 
We have
$$
z'^iaz'^{-i}\eq (1+b'z^q+c)z^iaz^{-i}(1+b'z^q+c)^{-1}\eq
(z^iaz^{-i}+b'z^qz^iaz^{-i}+cz^iaz^{-i})(1+b'z^q+c)^{-1}
\eq
$$
$$
({}^{(z)}\alpha^i (a)+{}^{(z,u)}_i\delta_{j}(a)z^j+\ldots +
{}^{(z,u)}_i\delta_{q+j}(a)z^{q+j}+\ldots +
b'z^q({}^{(z)}\alpha^i (a)+{}^{(z,u)}_i\delta_{j}(a)z^j+\ldots 
))(1+b'z^q+c)^{-1}\eq
$$
$$
({}^{(z)}\alpha^i (a)+b'{}^{(z)}\alpha^{q+i}(a)z^q+
{}^{(z)}\alpha^{i}(a)c+
{}^{(z,u)}_i\delta_{j}(a)z^j+\ldots
+{}^{(z,u)}_i\delta_{q+j}(a)z^{q+j}+\ldots +
$$
$$
b'\sum_{k\eq 
1}^{q}({}^{(z)}\alpha^{q-k}({}^{(z,u)}\delta_{j}({}^{(z)}\alpha^{k+i-1}(a) ) ))
z^{q+j}+
b'({}^{(z)}\alpha^{q}({}^{(z,u)}_i\delta_{j}(a)))
z^{q+j}
+\ldots )(1+b'z^q+c)^{-1}\eq
$$
$$
{}^{(z)}\alpha^i (a)+[{}^{(z,u)}_i\delta_{j}(a)z^j+\ldots +
{}^{(z,u)}_i\delta_{q+j}(a)z^{q+j}+\ldots +
b'\sum_{k\eq 
1}^{q}({}^{(z)}\alpha^{q-k}({}^{(z,u)}\delta_{j}({}^{(z)}\alpha^{k+i-1}(a) ) ))
z^{q+j}+
$$
$$
b'({}^{(z)}\alpha^{q}({}^{(z,u)}_i\delta_{j}(a)))
z^{q+j}
+\ldots )](1-b'z^q-c+\ldots )\eq
$$
$$
{}^{(z)}\alpha^i (a)+{}^{(z,u)}_i\delta_{j}(a)z^j+\ldots +
{}^{(z,u)}_i\delta_{q+j}(a)z^{q+j}+\ldots +
b'\sum_{k\eq 
1}^{q}({}^{(z)}\alpha^{q-k}({}^{(z,u)}\delta_{j}({}^{(z)}\alpha^{k+i-1}(a) ) ))
z^{q+j}+
$$
$$
b'({}^{(z)}\alpha^{q}({}^{(z,u)}_i\delta_{j}(a))
z^{q+j}
+\ldots
 -{}^{(z,u)}_i\delta_{j}(a){}^{(z)}\alpha^{j}(b')z^{q+j}+\ldots \eq
$$
$$
{}^{(z)}\alpha^i (a)+\ldots  +{}^{(z,u)}_i\delta_{q+j-1}(a)z'^{q+j-1}+
({}^{(z,u)}_i\delta_{q+j}(a)+b'{}^{(z)}\alpha^q({}^{(z,u)}_i\delta_{j}(a))-
{}^{(z,u)}_i\delta_{j}(a){}^{(z)}\alpha^j(b')
$$
$$
+b'\sum_{k\eq 
1}^{q}{}^{(z)}\alpha^{q-k}({}^{(z,u)}\delta_{j}({}^{(z)}\alpha^{k+i-1}(a))) -
{}^{(z,u)}_i\delta_{j}(a)\sum_{k\eq 0}^{j-1}{}^{(z)}\alpha^{k}(b))z'^{q+j} 
+\ldots ,
$$
since $z'^j\eq z^j+\sum_{k\eq 0}^{j-1}{}^{(z)}\alpha^{k}(b)z^{q+j}+\ldots $.

(iii)
We have
$$
z'az'^{-1}\eq bzaz^{-1}b^{-1}\eq
{}^{(z)}\alpha (a)+ b{}^{(z,u)}\delta_{j}(a){}^{(z)}\alpha^j(b^{-1})z^j+ 
\ldots\eq
$$
$$
{}^{(z)}\alpha (a)+ {}^{(z,u)}\delta_{j}(a){}^{(z)}\alpha (b^{-1})\ldots
{}^{(z)}\alpha^j(b^{-1})z'^j+ \ldots ,
$$
since ${}^{(z')}\alpha |_{Z(\bar{D})}={}^{(z)}\alpha |_{Z(\bar{D})}$.\\
$\Box$ 

\begin{corol}
\label{ozamene3}
In the situation of lemma \ref{ozamene} we have 
$$
j=w(xu(a)x^{-1}-u(a)),
$$
where $x\in D$ is any element with $w(x)=i$, if $a\in Z(\bar{D})$, $\alpha 
(a)=a$ and $(i,p)=1$, where $p=char D$. 

If $i=1$, we will denote $j$ by 
$j(u,a)$ or by $i(u,a)$.
\end{corol}

{\bf Proof.} Since for some parameter $z$ we have $x=b(1+x_1z+\ldots )z^i$, 
where $b, x_k\in u(\bar{D})$, the proof is easily follows from the proof of 
(ii) in lemma \ref{ozamene}.\\
$\Box$

In the sequel we will need the following definition.

\begin{defi}
Let
$(\alpha ,\beta )$ be endomorphisms of a division algebra $D$. A map
$\delta :$ $D\rightarrow D'$, where $D\subset D'$ are algebras, is called a
 $(\alpha ,\beta )$-derivation if it is linear and satisfy the following
identity

$$
\delta (ab)\eq \delta (a)\alpha (b)+\beta (a)\delta (b)
$$
where $a,b\in D$.\\
We will say that $(\alpha ,1)$-derivation  is an
$\alpha$-derivation.
\end{defi}

\begin{lemma}{(cf. \cite{Zh}, lemma 4)}
\label{lemma2}
Let $\delta$ be an $(\alpha ,\beta )$-derivation of an arbitrary division 
algebra $D$ such that
$\alpha ,\beta$ preserve $Z(D)$ and $\alpha|_{Z(D)}\ne \beta|_{Z(D)}$.

Then $\delta$ is an inner derivation, i.e. there exists $d\in D$ such
that
$$
\delta (a)\eq d\alpha (a)-\beta (a)d
$$
for all $a\in D$.
\end{lemma}

{\bf Proof.}
Put $d\eq \delta (a)(a^{\alpha}-a^{\beta})^{-1}$, where $a\in Z(D)$ is any 
element such
that $\alpha (a)\ne \beta (a)$. Put $\delta_{in}(x)\eq
d\alpha (x)-\beta (x)d$. We claim that $\delta\eq \delta_{in}$. Indeed,
consider the map $\bar{\delta}\eq \delta -\delta_{in}$. It is an
$(\alpha ,\beta )$-derivation. Take arbitrary $b\in D$. Then
$\bar{\delta}(ab)\eq \bar{\delta}(ba)$. But we have
$$\bar{\delta}(ab)\eq \bar{\delta}(a)\alpha (b)+\beta (a)\bar{\delta}(b)\eq
\beta (a)\bar{\delta}(b),$$
and
$$\bar{\delta}(ba)\eq \bar{\delta}(b)\alpha (a)+\beta (b)\bar{\delta}(a)\eq
\alpha (a)\bar{\delta}(b)$$
Therefore, $\bar{\delta}(b)\eq 0$ for any $b$.\\
$\Box$

\begin{prop}{(cf. \cite{Zh}, lemma 10)}
\label{X}
Let $D$ be a splittable division algebra. Let $n=Gal 
(Z(\bar{D})/\overline{Z(D)})$.
There exists a
parameter $z'$ such that
$$
{}^{(z',u)}_m\delta_j=0
$$
if $n\not | j$.
\end{prop}

{\bf Proof.} Since for $n=1$ there is nothing to prove, we will assume that 
$n>1$. Let $z$ be some fixed parameter. By \cite{JW}, prop. 1.7 
${}^{(z)}\alpha |_{Z(\bar{D})}$ has order $n$.

By proposition \ref{flyii}, ${}^{(z,u)}\delta_1$ is a 
$({}^{(z)}\alpha^2,{}^{(z)}\alpha
)$-derivation.
Since $n>1$, ${}^{(z)}\alpha^2|_{Z(\bar{D})}\ne {}^{(z)}\alpha 
|_{Z(\bar{D})}$. Therefore, by lemma \ref{lemma2},
${}^{(z,u)}\delta_1$ is an inner derivation and
${}^{(z,u)}\delta_1(a)\eq
 d{}^{(z)}\alpha^2(a)- {}^{(z)}\alpha (a)d$, $a\in \bar{D}$.
Put $z_1\eq  z-u(d)z^2$. By lemma \ref{ozamene}, (i) we have for any $a\in 
\bar{D}$
${}^{(z_1,u)}\delta_1(a)=0$ and ${}^{(z)}\alpha (a)={}^{(z_1)}\alpha (a)$. So,
${}^{(z_1,u)}\delta_1=0$ and ${}^{(z)}\alpha ={}^{(z_1)}\alpha$.

By proposition \ref{flyii}, ${}^{(z_1,u)}\delta_2$ is a 
$({}^{(z_1)}\alpha^3,{}^{(z_1)}\alpha )$-derivation.
 If $n\ne 2$ then it is inner and we can apply lemma
\ref{ozamene}.
 By induction we get that there exists a parameter $z_{n-1}$ such
that
${}^{(z_{n-1},u)}\delta_j=0$ for $j< n$ and ${}^{(z)}\alpha 
={}^{(z_{n-1})}\alpha$. It is easy to see that then 
${}^{(z_{n-1},u)}_m\delta_j=0$ for $j< n$ and all $m\in\dz$.
Note that ${}^{(z_{n-1},u)}\delta_n$ is a 
$({}^{(z_{n-1})}\alpha^{n+1},{}^{(z_{n-1})}\alpha )\eq
({}^{(z_{n-1})}\alpha ,{}^{(z_{n-1})}\alpha
)$-derivation, i.e. ${{}^{(z_{n-1},u)}\delta_{n}}{}^{(z_{n-1})}\alpha^{-1}$ is 
a derivation.

Note that  ${{}^{(z_{n-1},u)}\delta_{n+1}}$ is a
$({}^{(z_{n-1})}\alpha^2,{}^{(z_{n-1})}\alpha )$-derivation.
This follows  by proposition \ref{flyii}, since ${}^{(z_{n-1},u)}_m\delta_j=0$ 
for $j< n$ and all $m\in\dz$.
So, by lemma \ref{lemma2}, ${}^{(z_{n-1},u)}\delta_{n+1}$ is an inner 
derivation.
Using lemma \ref{ozamene}, (i) with $
z_{n+1}\eq z_{n-1}+bz_{n-1}^{n+2}$ for an appropriate $b$, we have
${}^{(z_{n+1},u)}\delta_j=0$ for $j< n+2$, $n\not |j$ and ${}^{(z)}\alpha 
={}^{(z_{n+1})}\alpha$. Moreover, ${}^{(z_{n+1},u)}_m\delta_j=0$ for $j< n+2$, 
$n\not |j$ and all $m\in\dz$. This easily follows from lemma \ref{triviall}.

By induction we can assume that there exists a
parameter $z_k$ such that
${}^{(z_{k},u)}_m\delta_j=0$ for $j< k+1$, $n\not |j$ and all $m\in\dz$, and 
${}^{(z)}\alpha ={}^{(z_{k})}\alpha$.

So, by proposition \ref{flyii}, if  $n\not | k+1$, then 
${{}^{(z_{k},u)}\delta_{k+1}}$ is an inner 
$({}^{(z_{k})}\alpha^{k+2},{}^{(z_{k})}\alpha )$-derivation. And
if $n | k+1$, we can apply the same arguments and conclude that
${}^{(z_{k},u)}\delta_{k+2}$
 is a $({}^{(z_{k})}\alpha^{k+2},{}^{(z_{k})}\alpha )$-derivation.
Therefore,  by lemma \ref{ozamene} there exists a parameter $z_{k+1}\eq
z_k+bz_k^{k+2}$
 ($z_k+bz_k^{k+3}$ if $n | k+1$) such that
${}^{(z_{k+1},u)}_m\delta_j=0$ for $j< k+2$, $n\not |j$ and all $m\in\dz$, and 
${}^{(z)}\alpha ={}^{(z_{k+1})}\alpha$ (or ${}^{(z_{k+1},u)}_m\delta_j=0$ for 
$j< k+3$, $n\not |j$ and all $m\in\dz$, and ${}^{(z)}\alpha 
={}^{(z_{k+1})}\alpha$ if $n | k+1$).

Since $z_{l+1}\eq
(1+b_lz_l^{k_l})z_l$ for every $l$, the sequence ${\{ z_l\}}_{l\eq 
1}^{\infty}$ converges
in $D$, which completes the proof of the proposition.\\
$\Box$

\begin{lemma}
\label{(5)}
Let $D$ be a splittable division algebra
as in proposition \ref{ooo}, of characteristic $p>0$.
Let $t\in Z(\bar{D})$ be an element such that $\alpha (t)=t$.

Let $j=i(u,t)$ be the minimal positive integer such that 
${}^{(z,u)}\delta_j|_{\sdf_p(t)}\ne
0$ (see corollary \ref{ozamene3}), and we assume $j<\infty$. Then the maps ${}^{(z,u)}_n\delta_m$, $kj\le 
m<(k+1)j$,
$k\in \{1,\ldots ,p-1 \}$ satisfy the following properties:

i) there exist elements $c_{n,m,k}\in\bar{D}$ such that
$$
{}^{(z,u)}_n\delta_m|_{\sdf_p(t)}=c_{n,m,1}\delta +\ldots +c_{n,m,k}\delta^k,
$$
where $\delta :\df_p(t)\rightarrow \df_p(t)$ is  a derivation such that 
$\delta (t)\eq 1$,  and 
$$c_{n,kj,k}\eq (k!)^{-1}{}^{(z,u)}_n\delta_j(t){}^{(z,u)}_{n+j}\delta_j(t)\ldots 
{}^{(z,u)}_{n+(k-1)j}\delta_j(t).
$$

ii) Let $\zeta =ord ({}^{(z)}\alpha |_{Z(\bar{D})})$. Then $\zeta |j$ and \\
$c_{n,kj,k}\ne 0$ if $(n, j )=1$ and 
${}^{(z)}\alpha ({}^{(z,u)}\delta_j(t))\ne {}^{(z,u)}\delta_j(t)$;\\
$c_{n,kj,k}\ne 0$ if ${}^{(z)}\alpha ({}^{(z,u)}\delta_j(t))= {}^{(z,u)}\delta_j(t)$ and 
$n , (n+j) , \ldots , (n+(k-1)j) \ne 0\mbox{\quad mod\quad}p$.

If ${}^{(z)}\alpha =id$, then $c_{n,kj,k}\ne 0$ iff 
$n , (n+j) , \ldots , (n+(k-1)j) \ne 0\mbox{\quad
mod\quad}p$.

\end{lemma}

{\bf Proof.} i) The proof is by induction on $k$. Let $a,b\in \df_p(t)$. For 
$k\eq 1$, by proposition \ref{flyii}, (ii) we have
$$
{}_n\delta_m(ab)\eq {}_n\delta_m(a)b+a{}_n\delta_m(b)
$$
because all the maps $\delta_q$, $q<j$ are equal to zero on $\df_p(t)$.
Hence, ${}_n\delta_m$ is a derivation on $\df_p(t)$, ${}_n\delta_m|_{\sdf_p(t)}=
c_{n,m,1}\delta$ and $c_{n,j,1}\eq
{}_n\delta_j(t)$.

For arbitrary $k$, by proposition \ref{flyii}, (i) and by the induction
hypothesis we have
$$
{}_n\delta_m(t^q)\eq q{}_n\delta_m(t)t^{q-1}+{}_n\delta_j(t)(\sum_{l\eq
0}^{q-2}(c_{n+j,m-j,1}\delta  +\ldots 
+c_{n+j,m-j,k-1}\delta^{k-1})(t^{q-1-l})t^l)+
$$
\begin{equation}
\label{(**)}
\ldots +{}_n\delta_{m-j}(t)(\sum_{l\eq  0}^{q-2}(c_{m-j+n,m-s,1}\delta 
)(t^{q-1-l})t^l).
\end{equation}
 Therefore, ${}_n\delta_m(t^p)\eq 0$,
because $k\le p-1$ and $\sum_{l\eq
0}^{p-2}\delta^i(t^{p-1-l})t^l\eq 0$ for $i\le p-2$. Hence,
${}_n\delta_m|_{\sdf_p(t)}\eq c_{n,m,1}\delta +\ldots +c_{n,m,p-1}\delta^{p-1}$ 
and we
only have to show that $c_{n,m,q}\eq 0$ for $q>k$.

Using (\ref{(**)}) we can calculate $c_{n,m,j}$. We have
$$
c_{n,m,1}\eq {}_n\delta_m(t);
$$
$$
c_{n,m,2}\eq \frac{1}{2!}({}_n\delta_m(t^2)-2c_{n,m,1}t)\eq \frac{1}{2}
({}_n\delta_j(t)(c_{n+j,m-j,1}\delta (t))+\ldots 
+{}_n\delta_s(t)(c_{s+n,m-s,1}\delta (t)))
$$
$$
\ldots
$$
$$
c_{n,m,q}\eq \frac{1}{q!}({}_n\delta_j(t)(\sum_{l\eq
0}^{q-2}c_{n+j, m-j, q-1}\delta^{q-1}(t^{q-1-l})t^l)+ \ldots
$$
$$
+{}_n\delta_{m-(q-1)j}(t)(\sum_{l\eq 0}^{q-2}c_{m+n-(q-1)j, (q-1)j, 
q-1}\delta^{q-1}
(t^{q-1-l})t^l))
$$
\begin{equation}
\label{recurrent}
\eq \frac{1}{q}({}_n\delta_j(t)c_{n+j, m-j, q-1}+ \ldots
+{}_n\delta_{m-(q-1)j}(t)c_{m+n-(q-1)j, (q-1)j, q-1})
\end{equation}
Hence, $c_{n, m, k+1}\eq \ldots \eq c_{n, m, p-1}\eq 0$ and 
$$c_{n, kj, k}\eq q^{-1}
{}_n\delta_j(t)c_{n+j, kj-j, k-1}\eq 
(k!)^{-1}{}^{(z,u)}_n\delta_j(t){}^{(z,u)}_{n+j}\delta_j(t)\ldots {}^{(z,u)}_{n+(k-1)j}\delta_j(t).
$$

ii) Let us prove first that $\zeta$ divide $i$. For, if $i$ is not divisible by $\zeta$, we have, by proposition \ref{flyii}, 
$$
{}^{(z,u)}\delta_j(tx)={}^{(z,u)}\delta_j(t){}^{(z)}\alpha^{j+1}(x)+{}^{(z)}\alpha 
(t)
{}^{(z,u)}\delta_j(x)={}^{(z,u)}\delta_j(xt)=
$$
$$
{}^{(z,u)}\delta_j(x){}^{(z)}\alpha^{j+1}
 (t) +{}^{(z)}\alpha (x){}^{(z,u)}\delta_j(t),
$$
where $x\in Z(\bar{D})$, $\alpha (x)\ne x$. But then 
${}^{(z)}\alpha^{j+1}(x)={}^{(z)}\alpha (x)$, a contradiction. 

If ${}^{(z)}\alpha =id$, the same arguments show that 
${}^{(z,u)}\delta_j(t)\in Z(\bar{D})$. 

If $x\in \bar{D}$ is an arbitrary element, this formulae shows ${}^{(z)}\alpha^{j}$ is an inner automorphism $ad({}^{(z,u)}\delta_j(t)^{-1})$. Therefore, 
${}^{(z)}\alpha^{j}({}^{(z,u)}\delta_j(t))={}^{(z,u)}\delta_j(t)$. 

Assume ${}^{(z)}\alpha ({}^{(z,u)}\delta_j(t))\ne {}^{(z,u)}\delta_j(t)$. It's clear then that 
$$
{}^{(z,u)}_{n+qj}\delta_j(t)=\sum_{l=0}^{n+qj-1}{}^{(z)}\alpha^l ({}^{(z,u)}\delta_j(t))\ne 0
$$ 
if $(n, j )=1$. So, $c_{n,kj,k}\ne 0$ by (i) in this case. 

If ${}^{(z)}\alpha ({}^{(z,u)}\delta_j(t))= {}^{(z,u)}\delta_j(t)$, then
${}^{(z,u)}_{n+qj}\delta_j(t)=(n+qj){}^{(z,u)}\delta_j(t)\ne 0$ iff $p$ does not divide 
$(n+qj)$. So, by (i) $c_{n,kj,k}\ne 0$ in this case iff $n , (n+j) , \ldots , (n+(k-1)j) \ne 0\mbox{\quad mod\quad}p$.

The lemma is proved.\\
$\Box$

\begin{lemma}
\label{ppp}
Let $D$ be a splittable division algebra as in lemma \ref{(5)}. Let $s\in 
Z(\bar{D})$ be an element such that $\alpha (s)=s$.
Let $i=i(u,s)$ be the minimal positive integer such that
${}^{(z,u)}\delta_i(s)\ne 0$ (see corollary \ref{ozamene3}).

If $p|i$, then for any positive integral $k$ there exists a map
${}^{(z,u)}\delta_{j(k)}$
such that ${}^{(z,u)}\delta_{j(k)}(s^{p^k})\ne
0$.
\end{lemma}

{\bf Proof.} We claim that ${}^{(z,u)}\delta_{p^qi}$ is the first map such that
${}^{(z,u)}\delta_{p^qi}|_{\sdf_p(s^{p^q})}\ne 0$. The proof is by induction on 
$q$. For
$q\eq 0$, there is nothing to prove. For arbitrary $q$, put $t\eq
s^{p^{q-1}}$. By proposition \ref{flyii} we have
$$
\delta_{p^qi}(t^p)\eq \delta_{p^{q-1}i}(t)\sum_{r\eq
0}^{p-2}{}_{1+p^{q-1}i}\delta_{p^{q-1}i(p-1)}(t^{p-1-r})t^r+ \sum_{l\eq
p^{q-1}i+1}^{p^qi-1}\delta_l(t)\sum_{r\eq
0}^{p-2}{}_{1+l}\delta_{p^qi-l}(t^{p-1-r})t^r
$$
 By induction and lemma \ref{(5)}, ${}_{1+l}\delta_{p^qi-l}|_{\sdf_p(t)}\eq
c_{1+l,p^qi-l,1}\delta +\ldots +c_{1+l,p^qi-1,p-2}\delta^{p-2}$ for 
$l>p^{q-1}i$. Therefore,
$\sum_{r\eq
0}^{p-2}{}_{1+l}\delta_{p^qi-l}(t^{p-1-r})t^r\eq 0$. By lemma \ref{(5)}, (ii),
${}_{1+p^{q-1}i}\delta_{p^{q-1}i(p-1)}|_{\sdf_p(t)}\eq
c_{1+p^{q-1}i, p^{q-1}i(p-1), 1}\delta +\ldots
+c_{1+p^{q-1}i, p^{q-1}i(p-1), p-1}\delta^{p-1}$ with
$c_{1+p^{q-1}i, p^{q-1}i(p-1), p-1}\ne 0$. Hence, $\delta_{p^qi}(t^p)\eq
-c_{1+p^{q-1}i, p^{q-1}i(p-1), p-1}\delta_{p^{q-1}i}(t)\ne 0$.

The same arguments show that ${}^{(z,u)}\delta_j(t^p)\eq 0$ for $j<p^qi$. So,
${}^{(z,u)}\delta_{p^qi}$ is the first non-zero map on $\df_p(s^{p^q})$. \\
$\Box$\\

\begin{lemma}
\label{svva}
Let $D$ be a splittable division algebra. Let $z$ be a fixed parameter and 
${}^{(z)}\alpha =id$, let $u$ be some fixed embedding 
$u:\bar{D}\hookrightarrow D$. 

Let ${}^{(z,u)}\delta_{i}$, $i\in \dn\cup\infty$ be the first non-zero map on 
$\bar{D}$. Assume $(i,p)=1$, where $p=char D$.  
Let ${}^{(z,u)}\delta_j$, $j>i$, $j\in \dn\cup\infty$ be the first map such 
that 
${}^{(z,u)}\delta_j\ne
0$ if $j$ is not divisible by $i$ and ${}^{(z,u)}\delta_j\ne
c_{j/i}{}^{(z,u)}\delta_i^{j/i}$  for some
$c_{j/i}\in \bar{D}$ otherwise. Then 

a) for $k< p=char D$ (arbitrary $k$ if $char D=0$) we have 
${}^{(z,u)}\delta_{ki}= 
c_{k}{}^{(z,u)}\delta_i^{k}$, where 
\begin{equation}
\label{(188)}
c_{k}\eq \frac{(i+1)\ldots (i(k-1)+1)}{k!},
\end{equation}
if $ki<j$.

b) if  condition (\ref{(188)}) is satisfied for any $k$ with $ki<j$, then 
${}^{(z,u)}_{-i}\delta_{q}=0$ for $i<q<j$ and ${}^{(z,u)}_{-i}\delta_{j}$ is a 
derivation.  

\end{lemma}

{\bf Remark.} We will call the number $i(u,z)=\min_{a\in 
\bar{D}}\{w(zu(a)z^{-1}-u(a))\}$ defined in this lemma {\it a local height}. 
The number $i=i(z,u)$ in lemma coinside with the level of $D$ defined in 
\cite{Sa} if $D$ has index $p=char D$ and $D$ is splittable. As it follows 
from lemmas \ref{ozamene}, \ref{ozamene2} (see below), $i(z,u)$ does not 
depend on $z,u$ in this case. Corollary \ref{ozamene3} completes then the 
proof that it coinside with the level defined by Saltman in the case $D$ is 
splittable.  This number will play an important role in this work. It was one 
of the important  parameters in \cite{Zh}. Recall the definition of {\it 
level}: $h(D)=\min \{w(ab-ba)-w(a)-w(b)\}$. 

{\bf Proof.} If we compare coefficients in formulae for $\delta_{ki}(ab)$ from 
proposition \ref{flyii} with coefficients in formulae for $\delta_i^k(ab)$ 
multiplied by $c_k$, we must have 
$$
c_kk=((k-1)i+1)c_{k-1},
$$ 
where from follows a). 

>From the other hand side, if ${}_{-i}\delta_q$, $q>i$ is the first nonzero 
map after 
${}_{-i}\delta_i$, it must be a derivation by proposition \ref{flyii}, (i). 
Note that in characterictic zero case this can happens only if $q\ge j$, because 
a map $c\delta_i^k$ can not be a derivation if $k>1$, which proves b) in 
this case. 

Since the maps $\delta_q$ are uniquely defined, by lemma \ref{triviall}, by 
the maps ${\tilde{\delta}}_l$, $l\le q$, and the maps ${\tilde{\delta}}_q$ are 
uniquely defined by the maps 
${}_{-i}{{\delta}}_l$, $l\le q$, and ${}_{-i}{{\delta}}_q$ are  linear 
combinations  
of ${{\delta}}_l$, $l\le q$ with integer coefficients, we see that b) holds in 
arbitrary characteristic.\\
$\Box$

{\bf Remark.} So we see that the maps ${}_i\delta_q$ in this lemma satisfy the 
same identities as $\delta_{q/i}$. This can be thought of as a possible 
reduction from level $i$ to level $1$.

\begin{defi}
\label{lemdef}
Let $D$ be a splittable division algebra. Let $u$ be some fixed embedding 
$u:\bar{D}\hookrightarrow D$. 
Let $s\in Z(\bar{D})$ be an element such that $\alpha (s)=s$.
Let $i=i(u,s)$ be the minimal positive integer such that
${}^{(z,u)}\delta_i(s)\ne 0$ ( corollary \ref{ozamene3} shows that $i$ does 
not depend on 
$z$). Assume $(i,p)=1$, where $p=char D$. 
Define 
$$
d(u,s)=\max_{z}\{w(z^{-i}u(s)z^i-u(s)-u({}^{(z,u)}_{-i}\delta_i(s))z^i)\} \in 
\dn\cup\infty ,
$$
\end{defi}

As we can see from lemma \ref{svva} b), $d(u,s)$ can be interpreted under some conditions as the number $j$ there. So, this definition was motivated by this lemma.

\begin{lemma}
\label{vtorinv}
In the definition above for $p=char D>0$ and 
${}^{(z)}\alpha |_{Z(\bar{D})}=id$ we have

i) ${d(u,s)}=2i \mbox{\quad mod\quad} p$ if $d(u,s)<\infty$;

ii) If ${}^{(z,u)}_{-i}\delta_i(s)\ne 0$, the map 
${}^{(z,u)}_{-i}\delta_{d(u,s)+(p-1)i}$ is the first map such that 
${}^{(z,u)}_{-i}\delta_{d(u,s)+(p-1)i}(s^p)\ne 0$ for any parameter $z$. In 
particular, if $d(u,s)=\infty$, $[u(s^p), z^i]=0$.

\end{lemma} 

{\bf Proof.}
(ii) Let ${}^{(z,u)}_{-i}\delta_{\kappa}$ be the first map such that 
${}^{(z,u)}_{-i}\delta_{\kappa}(s^p)\ne 0$. By corollary \ref{ozamene3}   
$\kappa$ does not depend on $z$. By the same reason, 
${}^{(z,u)}_{-i}\delta_i$ is the first map such that 
${}^{(z,u)}_{-i}\delta_i(s)\ne 0$ for any $z$. 

Put $w:\eq d(u,s)+(p-1)i$ and fix $u,z$. By proposition \ref{flyii} we have 
$$
{}_{-i}{{\delta}}_w(s^p)\eq 
{}_{-i}{{\delta}}_{d(u,s)}(s)\sum_{q\eq 0}^{p-2}
{}_{d(u,s)-i}{{\delta}}_{(p-1)i}(s^{p-1-q})s^q+
$$
$$
\sum_{k\eq d(u,s)+1}^{w-1} {}_{-i}{{\delta}}_k(s)\sum_{q\eq 0}^{p-2}
{}_{k-i}{{\delta}}_{w-k}(s^{p-1-q})s^q
$$ 
 By lemma \ref{(5)}, ${}_{k-i}{{\delta}}_{w-k}|_{\sdf_p(s)}\eq 
c_{k-i,w-k,1}\delta
+\ldots +c_{k-i,w-k,p-2}\delta^{p-2}$ for $w-k<(p-1)i$ and 
${}_{d(u,s)-i}{{\delta}}_{(p-1)i}|_{\sdf_p(s)}\eq c_{d(u,s)-i,(p-1)i,1}\delta 
+\ldots
+c_{d(u,s)-i,(p-1)i,p-1}\delta^{p-1}$ with $c_{d(u,s)-i,(p-1)i,p-1}\ne 0$ if $d(u,s)-i=i \mbox{\quad mod\quad} p$. Indeed, as we have shown in the proof of lemma \ref{(5)}, (ii), 
the order $n$ of the automorphism ${}^{(z)}\alpha$ on ${}^{(z,u)}\delta_i(s)$  must divide $i$, so $(n,p)=1$. Now we have two possibilities: $n{\not |}d(u,s)$ and $n|d(u,s)$. 

In the first case we can repeat the arguments to the first assertion in lemma \ref{(5)}, (ii) to show that $c_{d(u,s)-i,(p-1)i,p-1}\ne 0$. In the second case we have 
${}_{d(u,s)-i+qi}\delta_i(s)= (d(u,s)-i+qi)/i{}_i\delta_i(s)\ne 0$ if $d(u,s)-i+qi$ is not divided by $p$. So, by lemma \ref{(5)}, (i) $c_{d(u,s)-i,(p-1)i,p-1}\ne 0$ iff $d(u,s)-i=i \mbox{\quad mod\quad} p$ in this case. 

Hence,
$$
{}_{-i}{{\delta}}_w(s^p)\eq
-{}_{-i}{{\delta}}_{d(u,s)}(s)c_{d(u,s)-i,(p-1)i,p-1}\ne 0
$$
if $d(u,s)-i=i \mbox{\quad mod\quad} p$.

This also shows that ${}_{-i}{{\delta}}_w$ is the {\it first} map
such that ${}_{-i}{{\delta}}_w|_{\sdf_p(s^p)}\ne 0$ if $d(u,s)-i=i \mbox{\quad 
mod\quad} p$.

i)  By Skolem-Noether theorem there exists a parameter $z'$ in $D$ such that 
${}^{(z')}\alpha =id$. Put 
$$d'(u,z',s)= w(z'^{-j}u(s)z'^{j}-u(s)-u({}^{(z',u)}_{-i}\delta_i(s))z'^i).$$
Since ${}^{(z')}\alpha =id$, the map 
${}^{(z',u)}\delta_i$ is the first map such that 
${}^{(z',u)}_{-i}\delta_i(s)\ne 0$. If 
$d'(u,z',s)\ne 2i \mbox{\quad mod\quad} p$, we can find a parameter $z''$ such 
that 
$d'(u,z'',s)>d'(u,z',s)$ using lemma \ref{ozamene}, (ii). Continuing this 
procedure, we find a parameter $z$ such that $d'(u,z,s)=2i \mbox{\quad 
mod\quad} p$ or $d'(u,z,s)=\infty$. 

Using arguments from ii) we get that the map 
${}^{(z,u)}_{-i}\delta_{d'(u,z,s)+(p-1)i}$ is the first map such that 
${}^{(z,u)}_{-i}\delta_{d'(u,z,s)+(p-1)i}(s^p)\ne 0$ for the parameter $z$. As 
it was noted in the beginning of the proof, the number $\kappa =d'(u,z,s)+(p-1)i$ does not depend on the parameter. 
Since $d'(u,z,s)\le d(u,s)$, we get $d'(u,z,s)= d(u,s)$. For, otherwise we can repeat the arguments from (ii) and conclude that ${}^{(z,u)}_{-i}\delta_{d(u,s)+(p-1)i}(s^p)= 0$, a contradiction.
The lemma is proved.\\
$\Box$

 It would be interesting to know more about a behaviour of 
${}^{(z,u)}_m\delta_j$ with respect to the embedding $u$. We will give an 
answer in one special case, namely,  when $\bar{D}=Z(\bar{D})$ and 
$Z(\bar{D})/\overline{Z(D)}$ is a simple extension.

\begin{lemma}
\label{simple}
Let $D$ be a division algebra such that $char D=p>0$,
$\bar{D}=Z(\bar{D})$, $Z(\bar{D})$ is not perfect and 
$Z(\bar{D})/\overline{Z(D)}$ is a simple extension (so, $D$ is splittable).
Let $\bar{u}$ be a primitive element of the extension 
$Z(\bar{D})/\overline{Z(D)}$ such that $\bar{u}\notin (Z(\bar{D}))^p$ and let 
$u$ be any lift of  $\bar{u}$ in $D$.

Then there exists an embedding
$u:\bar{D}\hookrightarrow D$ such that $u(\bar{u})=u$ and any map 
${}^{(z,u)}_m\delta_j$ is uniqely defined by the values 
${}^{(z,u)}_m\delta_j(u^q)$ or, equivalently, by the values 
${}^{(z,u)}_l\delta_k(u)$, $k\le j$.

In particular, if ${}^{(z,u)}_m\delta_k(u)=0$ for $k\le j$, then 
${}^{(z,u)}_m\delta_j=0$.
\end{lemma}

{\bf Proof.} Consider a field $Z(D)(u)$. It is a complete discrete valued 
field as a finite extension of $Z(D)$. By classical Cohen theorem, there 
exists an embedding $\overline{Z(D)(u)}=\bar{D}\hookrightarrow Z(D)(u)\subset 
D$.
By \cite{Co}, lemmas 11,12 the embedding is completely defined by a $p$-basis 
$\Gamma$ of the field $\overline{Z(D)(u)}$. Namely, for any lift $G$ of a 
given $p$-basis $\Gamma$ there exists an embedding $s$ such that $G\subset 
s(\overline{Z(D)(u)})$.

Let's show that there exists a $p$-basis $\Gamma$ of the field $\bar{D}$ such 
that $\bar{u}\in \Gamma$ and $\Gamma\ni \gamma\in \overline{Z(D)}$ if 
$\gamma\ne \bar{u}$.

Consider a set of all non-void sets $\Gamma'$ of elements $\gamma_{\tau}\in 
\bar{D}$ satisfying the following property:\\
A) $\bar{u}\in\Gamma'$, $\Gamma' \ni \gamma\in \overline{Z(D)}$ if $\gamma\ne 
\bar{u}$ and $[{\bar{D}}^p(\gamma_1,\ldots ,\gamma_r):\bar{D}^p]=p^r$ for any 
$r$ distinct elements of $\Gamma'$.

This set is not void, since it contains the set $\Gamma'=\{\bar{u}\}$. By 
Zorn's lemma,
there exists a maximal set $\Gamma$ satisfying A). Then 
$\bar{D}=\bar{D}^p(\Gamma )$. Indeed, since $\overline{Z(D)}^p(\bar{u})\subset 
\bar{D}^p(\Gamma )$, it suffice to show that any element from 
$\overline{Z(D)}$ lies in $\bar{D}^p(\Gamma )$. Suppose $a\in 
\overline{Z(D)}$, $a\notin \bar{D}^p(\Gamma )$. Then  the set $\Gamma'=\{a\cup 
\Gamma\}$ satisfy A), a contradiction with maximality of $\Gamma$.

Now, we can take a lift of $\Gamma$ in the following way. We take $u$ as a 
lift of $\bar{u}$, and we take lifts of all other elements in $Z(D)$. This 
lift defines an embedding $u:\bar{D}\hookrightarrow D$.

Let us show that any map ${}^{(z,u)}_m\delta_j$ (for some fixed $z$) is 
uniqely defined by the values ${}^{(z,u)}_l\delta_k(u)$, $k\le j$. We have 
$u(\bar{D})=u(\overline{Z(D)})(u)$ and any element $a\in u(\bar{D})$ can be 
represented as a polynomial in finite number of elements from $\Gamma$ with 
coefficients from $u(\bar{D})^{p^k}$ for any $k>0$.

Note that for any $j$ there exists $k>0$ such that for any $b\in 
\overline{Z(D)}^{p^k}$ ${}^{(z,u)}_l\delta_q(b)=0$ for all $q\le j$ and all 
$l$. Indeed, assume  ${}^{(z,u)}_1\delta_q(b)\ne 0$ for some  $q\le j$, $b\in 
\overline{Z(D)}^{p^k}$ and  ${}^{(z,u)}_l\delta_s(c)= 0$ for all $l$, all 
$c\in \overline{Z(D)}^{p^k}$ and all $s<q$. Then, since
${}^{(z)}\alpha |_{\overline{Z(D)}}=id$ and by proposition \ref{flyii},
${}^{(z,u)}_l\delta_s(b^p)= 0$ for all $b\in \overline{Z(D)}^{p^k}$, all $l$ 
and all $s\le q$.

Now, since $u(\bar{D})^{p^k}=u(\overline{Z(D)})^{p^k}(u^{p^k})$, any element 
$a\in u(\bar{D})$ can be represented as a polynomial in finite number of 
elements from $\Gamma$ with coefficients from $u(\overline{Z(D)})^{p^k}$. 
Since all elements except $u$ in $\Gamma$ belong to the center $Z(D)$, the 
value of ${}^{(z,u)}_m\delta_j(a)$
is uniqely determined by the values ${}^{(z,u)}_m\delta_j(u^l)$ that are 
uniqely defined, by proposition \ref{flyii}, by the values 
${}^{(z,u)}_l\delta_k(u)$, $k\le j$. \\
$\Box$

{\bf Remark} In the case $Z(\bar{D})$ perfect field there is only one 
embedding $u$, which is compatible with the embedding 
$\overline{Z(D)}\hookrightarrow Z(D)$. So, the assertion of lemma is easy in 
this case.

\begin{lemma}{(cf. \cite{Zh}, lemma 8)}
\label{ozamene2}

In the situation of lemma \ref{simple}
suppose
${}^{(z,u)}_m\delta_1\eq\ldots \eq {}^{(z,u)}_m\delta_{j-1}\eq 0$, 
${}^{(z,u)}_m\delta_j\ne 0$. Let $n$ be the order of ${}^{(z)}\alpha$. Then

(i) for $u'\eq u+bz^q$, $b\in u(\bar{D})$, $n|q$ we have
${}^{(z,u')}_m\delta_l={}^{(z,u)}_m\delta_l$, $l<q$ and
$${}^{(z,u')}_m\delta_q (\bar{u})\eq {}^{(z,u)}_m\delta_q 
(\bar{u})+{}^{(z)}\alpha^m
(\bar{b})
-\frac{\partial}{\partial \bar{u}}({}^{(z)}\alpha^m (\bar{u}))\bar{b},$$
where the derivative is taken in the field $\bar{D}=\bar{D}^p(\Gamma )$.

(ii) Suppose ${}^{(z)}\alpha \eq id$. Then for $u'\eq u+bz^q$, $b\in 
u(\bar{D})$  we have
${}^{(z,u')}_m\delta_l={}^{(z,u)}_m\delta_l$, $l<q+j$ and
$${}^{(z,u')}_m\delta_{q+j}(\bar{u})\eq {}^{(z,u)}_m\delta_{q+j}(\bar{u})+
{}^{(z,u)}_m\delta_{j}(\bar{b})
-\frac{\partial }{\partial \bar{u}}({}^{(z,u)}_m\delta_j (\bar{u}))\bar{b}, $$
where the derivative is taken in the field $\bar{D}=\bar{D}^p(\Gamma )$.

(iii) Suppose ${}^{(z)}\alpha\eq id$. Let $\bar{u'}\in \bar{D}$ be any 
primitive element of the extension $\bar{D}/\overline{Z(D)}$ satisfying the 
conditions of lemma
\ref{simple}, and let $u'\in D$ be any lift of $\bar{u'}$.
Then  we have
${}^{(z,u')}_m\delta_l={}^{(z,u)}_m\delta_l$, $l<j$ and
$$
{}^{(z,u')}_m\delta_{j}(\bar{u'})\eq {}^{(z,u)}_m\delta_{j}(\bar{u})
\frac{\partial }{\partial \bar{u}}(\bar{u'}),
$$
where the derivative is taken in the field $\bar{D}=\bar{D}^p(\Gamma )$.
\end{lemma}

{\bf Proof.}
First of all, let's note that there exists $k\in \dn$ such that for any $a\in 
\overline{Z(D)}^{p^k}$ holds $u(a)-u'(a)=0 \mbox{\quad mod\quad} M_D^{q+1}$, 
where $u'$ is any another embedding, $q\in \dn$ is any given number.

Indeed, assume for any $c\in \overline{Z(D)}^{p^s}$ holds $u(c)-u'(c)=0 
\mbox{\quad mod\quad} M_D^{l}$, i.e.
$u(c)=u'(c)+c_lz^l+\ldots$, where  $c_l\in u'(\bar{D})$.  Then 
$u(c^p)=(u(c))^p=
(u'(c))^p+pu'(c)^{p-1}c_lz^l+\ldots$, so
$u(c^p)-u'(c^p)=0 \mbox{\quad mod\quad} M_D^{l+1}$.

>From this immediately follows that $u(a)-u'(a)=0 \mbox{\quad mod\quad} 
M_D^{q}$ for any $a\in \bar{D}$ if $u'$ is defined by the element $u'=u+bz^q$, 
because $u(\bar{u})-u'(\bar{u})=bz^q$. Moreover, if we represent $a$ as some 
polynomial $P(\gamma_1,\ldots ,\gamma_r, \bar{u})$ with coefficients from
$\overline{Z(D)}^{p^k}$, then it is clear that
$$[u(a)-u'(a)]z^{-q}=-\frac{\partial}{\partial \bar{u}}(P(\gamma_1,\ldots 
,\gamma_r, \bar{u}))
\bar{b} \mbox{\quad mod\quad} M_D
$$
if $n|q$,
since $u(\gamma_l)=u'(\gamma_l)$ for any $l$ and $z^quz^{-q}=u \mbox{\quad 
mod\quad} M_D$. It is also clear  that the derivative can be taken even in the 
field $\bar{D}^p(\Gamma )$.
So, we have \\
(i)
$$
z^mu'z^{-m}\eq z^m(u+bz^q)z^{-m}\eq
u({}^{(z)}\alpha^m (\bar{u}))+u({}^{(z,u)}_m\delta_j(\bar{u}))z^j +\ldots +
(u({}^{(z)}\alpha^m(\bar{b}))
$$
$$
+u({}^{(z,u)}_m\delta_j(\bar{b}))z^j+ \ldots )z^q\eq
u({}^{(z)}\alpha^m(\bar{u}))+\ldots +(u({}^{(z,u)}\delta_q (\bar{u}))+
u({}^{(z)}\alpha^m(\bar{b})))z^q+\ldots \eq
$$
$$
u'({}^{(z)}\alpha^m (\bar{u}))+\ldots +(u'({}^{(z,u)}_m\delta_q 
(\bar{u}))+u'({}^{(z)}\alpha^m(\bar{b}))-
u'(\frac{\partial}{\partial \bar{u}}{}^{(z)}\alpha^m 
(\bar{u})\bar{b}))z^q+\ldots ,
$$
(ii)
We have
$$
z^mu'z^{-m}\eq z^m(u+bz^q)z^{-m}\eq
u(\bar{u})+u({}^{(z,u)}_m\delta_j (\bar{u}))z^j+\ldots +(u(\bar{b})+
u({}^{(z,u)}_m\delta_j(\bar{b}))z^j+ \ldots )z^q\eq
$$
$$
u(\bar{u})+u({}^{(z,u)}_m\delta_j (\bar{u}))z^j+\ldots
+(u({}^{(z,u)}_m\delta_q 
(\bar{u}))+u(\bar{b}))z^q+u({}^{(z,u)}_m\delta_{q+1}(\bar{u}))
z^{q+1}+\ldots
$$
$$
+u({}^{(z,u)}_m\delta_{q+j-1}(\bar{u}))z^{q+j-1}+(u({}^{(z,u)}_m\delta_{q+j}(\bar{u}))+
u({}^{(z,u)}_m\delta_j (\bar{b})))z^{q+j}+\ldots\eq
$$
$$
u'(\bar{u})+u'({}^{(z,u)}_m\delta_j (\bar{u}))z^j+\ldots 
+u'({}^{(z,u)}_m\delta_{q+j-1}(\bar{u}))z^{q+j-1}+
(u'({}^{(z,u)}_m\delta_{q+j}(\bar{u}))+u'({}^{(z,u)}_m\delta_j (\bar{b}))-
$$
$$
u'(\frac{\partial}{\partial \bar{u}}({}^{(z,u)}_m\delta_j 
(\bar{u}))\bar{b}))z^{q+j}+\ldots
$$
(iii) Assume $u'=u(\bar{u'})+a_1z+\ldots$, where $a_i\in u(\bar{D})$. Since, 
by proposition \ref{flyii}, the map  ${}^{(z,u)}_m\delta_j$ is a derivation,
we have
$$
z^mu'z^{-m}\eq [u(\bar{u'})+u({}^{(z,u)}_m\delta_j(\bar{u'}))z^j+\ldots ] +
[a_1+u({}^{(z,u)}_m\delta_j(a_1)z^j+\ldots ]z+\ldots =
$$
$$
u'+u({}^{(z,u)}_m\delta_j(\bar{u'}))z^j+\ldots 
=u'+u({}^{(z,u)}_m\delta_j(\bar{u})
\frac{\partial }{\partial \bar{u}}(\bar{u'}))z^j+\ldots =
u'+u'({}^{(z,u)}_m\delta_j(\bar{u})
\frac{\partial }{\partial \bar{u}}(\bar{u'}))z^j+\ldots
$$
$\Box$

\section{The period-index problem}

In this section we will prove the following theorem.

\begin{th}
\label{gipoteza}
The following conjecture: the exponent of $A$ is equal to its index for any 
division algebra $A$  over a $C_2$-field $F$  has the positive answer for 
$F\eq F_1((t))$, where $F_1$ is a  $C_1$-field.  
\end{th}

Recall that a field $F$ is called a {\it $C_i$-field } if any homogeneous form 
$f(x_1, \ldots ,x_n)$ 
of degree $d$ in $n>d^i$ variables with coefficients in $F$ has a non-trivial 
zero. 
Some basic properties of $C_i$-fields see, for example, in \cite{PY}.

This conjecture was proposed by M. Artin and was solved for some another 
examples of the field $F$ by many authors. As it is known for me, the positive 
answer  for all division algebras of index $ind A= 2^a3^b$ was given in  
\cite{PY}, for  division algebras over the field $F=k((X))((Y))$, where $k$ is 
a perfect field of characteristic $p\ne 0$ such that $\dim_{\sdf_p}k/\wp 
(k)=1$, was given by Tignol in the Appendix in \cite{AJ} (we include this case though $F$ may not be a $C_2$-field), for division 
algebras of index prime to the characterictic of $F$, where $F$ is a function 
field of a surface, was given in \cite{DJ}. I propose, the positive answer was 
also known for division algebras over $F=F_1((t))$ of characteristic 0.  We 
will give the prove of the theorem above in any characteristic.

{\bf Proof.} 1) Recall that any extension of a $C_1$-field is simple. 
Indeed, suppose $E\eq \bar{F}(u_1, \ldots , u_r)$. Consider the field $K\eq 
\bar{F}(u_1^p, \ldots , u_r^p)$.
By Tsen's theorem, $K$ and $E$ are $C_1$-fields. So, the form 
$x_1^p+x_2^pu_1+\ldots + 
x_p^pu_1^{p-1}+x_{p+1}^pu_2$ has a non-trivial zero in $E$. But $x_i^p\in K$ 
and elements 
$1, u_1,\ldots , u_1^{p-1}, u_2$ are linearly independent over $K$, a 
contradiction.

2) Assume the theorem is known in the prime exponent case. We deduce the 
theorem  by ascending induction on $e=exp A$. If $e$ is not a prime number, 
then write $e=lm$. By assumption $A^{\otimes m}$ can be split by a field 
extension $F\subset F'$ of degree $l$. This implies that $A_{F'}$  has 
exponent dividing $m$. Note that $F'$ is also a Laurent series field. By the 
induction hypothesis applied to the pair $(F', A_{F'})$, there exists a field 
extension $F'\subset L$ of degree dividing $m$ splitting $A_{F'}$. Therefore 
$A$ is split by the extension $F\subset L$ of degree dividing $lm$ and we 
conclude the theorem. 

3) So, let $exp A= l$ be a prime number. By the basic properties of the 
exponent and the index (see, e.g. \cite{PY}) we have then $ind A= l^k$ for 
some natural $k$. 

Suppose $(l,p=char F)=1$.

It is known that the conjecture is true for all division algebras of index 
$ind A= 2^a3^b$, so we can assume $l\ne 2,3$. We can assume  $F$ contains the 
group $\mu_l$ of $l$-roots of unity, because 
$[F(\mu_l):F]<l$ and we can reduce the problem to the algebra $A\otimes_F 
F(\mu_l)$. Then by the Merkuriev-Suslin theorem 
$A$ is similar to the tensor product of symbol-algebras of index $l$. 

 To conclude the statement of the corollary it is sufficient to prove that 
every two symbol algebras $A_1, A_2$ contain $F$-isomorphic maximal subfields. 

Since every division algebra over a $C_1$-field is trivial and every field 
extension is simple, every symbol-algebra of index $l$ over $F$ is  
splittable. Since $(l,p)=1$, it is good splittable  and its residue field is a 
cyclic Galois extension of $\bar F$.
So, if $z_i$ is a parameter from proposition \ref{X} for algebra $A_i$, then 
$z_i$ acts on $\bar{A_i}$ as a Galois automorphism and $z_i^l\in F$. We have 
$v(z_i^l)=1$ ($v$ is the valuation on $F$).

Let us show that $A_1$ contains a $l$-root of any element $u$ in $F$ with 
$v(u)\ne 0$. So, $A_1$ will contain a subfield isomorphic to $F(z_2)$. Since 
for any element $1+b$, $v(b)>0$  there exists 
a $l$-root $(1+b)^{1/l}\in F$, it is sufficient to prove that $A_1$ contains 
any $l$-root of elements $ct$, $c\in u(\bar{F})$, where $u$ is some fixed 
embedding $u:\bar{A_1}\hookrightarrow A_1$. 

Assume  $z_1^l=c_1t$, $c_1\in u(\bar{F})$. Note that for any element $b\in 
u(\bar{A_1})$ we have $(bz_1)^l=u(N_{\bar{A_1}/\bar{F}}(b))z_1^l$. But the 
norm map $N_{\bar{A_1}/\bar{F}}$ is surjective, since $\bar F$ is a 
$C_1$-field (see, e.g. \cite{PY},  3.4.2), so there exists $b$ such that 
$(bz)^l=ct$.

4) Suppose now $exp A=p$. Then $ind A=p^k$. 

By Albert's theorem (in \cite{Al}) there exists a field $F'=F(u_1^{1/p},\ldots 
,u_k^{1/p})$ which splits $A$. Using the same arguments as in 1) one can show 
that every such a field has maximum two generators, say $F'=F(u_1^{1/p}, 
u_2^{1/p})$. Therefore, $ind A\le p^2$. If $ind A=p$, there is nothing to 
prove, so we assume $ind A=p^2$ and $F'$ is a maximal subfield in $A$. 

5) Suppose $F_1$ is a perfect field. 

By Albert's theorem, $A\cong A_1\otimes_F A_2$, where $A_1,A_2$ are cyclic 
algebras of degree $p$, $A_1=(L_1/F,\sigma_1, u_1)$, $A_2=(L_2/F, \sigma_2, 
u_2)$. Since $F_1$ is perfect, $\bar{A_1}/\bar{F}$,
$\bar{A_2}/\bar{F}$ are Galois extensions. So, $A_1, A_2$ are good splittable. 
Let us show that $A_1, A_2$ have common splitting field of degree $p$ over 
$F$. 
This leads to a contradiction.

By proposition \ref{X} there exist parameters $z_1\in A_1$, $z_2\in A_2$ such that 
they act on $\bar{A_1}$, $\bar{A_2}$ as Galois automorphisms. Note that then 
$z_1^p, z_2^p\in F$. Let us show that $F(z_1)$ splits $A_2$.

Consider the centralizer $D=C_A(F(z_1))$. Consider the element 
$t_1=z_2z_1^{-1}$. We have $t_1^p\in F$, $w(t_1)=0$, where 
$w$ denote the unique extension of the valuation $v$ on $F$. Since 
$\bar{D}/\overline{Z(D)}$ is a Galois extension, there exists an element 
$b_1\in F$ such that $w(t_1-b_1)>0$. Since $(t_1-b_1)^p\in F$, there exists 
natural $k_1$ such that 
$w((t_1-b_1)z_1^{-k_1})=0$. Denote $t_2=(t_1-b_1)z_1^{-k_1}$. We have again 
$t_2^p\in F$. Repeating this arguments and using the completeness of $D\subset 
A$ we get \\
$z_2=t_1z_1=(t_2z_1^{k_1}+b_1)z_1=\ldots =b_1z_1+b_2z_1^{k_1+1}+\ldots$,\\
so, $z_2\in F(z_1)=Z(D)$. 
 
6) Suppose $F_1$ is not perfect. 

Since $F'$ is generated by two elements over $F$, it contains all $p$-roots of 
$F$. Then, every two elements $u,z\in F$ such that $z^{1/p}\notin F(u^{1/p})$, 
where 
$z^{1/p}, u^{1/p}\in F'$, also generate $F'$ over $F$. This follows from the 
same arguments as in 1), 4). 
 
Now take $u\in F_1\backslash F_1^p$, $z=u+t$. It's clear that $p$-roots of 
these elements generate $F'$ over $F$. Moreover, the fields $F(u^{1/p}), 
F(z^{1/p})$ are {\it "unramified"} over $F$, i.e. 
$[\overline{F(u^{1/p})}:\bar{F}]=p=[F(u^{1/p}):F]$, 
$[\overline{F(z^{1/p})}:\bar{F}]=p$.  Denote $u_1=u^{1/p}$,
$u_2=z^{1/p}$ in $F'$. Then by Albert's theorem, $A\cong A_1\otimes_F A_2$, 
where $A_1,A_2$ are cyclic algebras of degree $p$, $A_1=(L_1/F,\sigma_1, u)$, 
$A_2=(L_2/F, \sigma_2, z)$. 

Concider the centralizer $D=C_A(F(u_1))$. Suppose $\bar{D}/\overline{Z(D)}$ is 
a separable extension. Then there exist a lift $u:\bar{D}\hookrightarrow D$ of 
arbitrary embedding $u':\overline{F(u_1)}\hookrightarrow F(u_1)$. Consider the 
embedding $u'=u_1$ defined in lemma \ref{simple}. 
Since $F(u_1)/F$ is a purely inseparable extension, $u'$ is a good embedding, 
so $u$ is a good embedding of $\bar{D}=\bar{A}$ in $D\subset A$. So, we get 
$A$ is a good splittable algebra, and $u(\bar{A})$ contain a purely 
inseparable over $F$ element.  But this is a contradiction with lemma 
\ref{ppp}. So, $\bar{A}/\bar{F}$ can not contain a separable subextension, 
because in this case  $\bar{D}/\overline{Z(D)}$ must be a separable extension.

Now we can use, for shorteness, lemmas A.4., A.6. of Tignol in Appendix to the 
paper 
\cite{AJ}. These lemmas show that a tensor product $A_1\otimes A_2$ of any two 
symbols 
$A_1, A_2$ is similar either to a single symbol in $Br(F)$ (in which case we are done) or to a  
product of two symbols of level zero. Recall that, by Saltman's results in 
\cite{Sa}, every  division 
algebra of level zero is tame, which means in our case that the residue 
division algebra is a separable extension over $\bar{F}$. A notion of level 
was already discussed above in remark to lemma \ref{svva}. 

So, assume $A\sim D_1\otimes D_2$, where $D_1, D_2$ are tame division algebras of 
degree $p$ over $F$. We can assume $A$ and $D_1\otimes D_2$ are  division algebras, so
$A\cong D_1\otimes D_2$.  Since $D_1, D_2$ are tame, 
we conclude $\bar{A}$  must contain a separable element, a contradiction.

The theorem is proved. \\
$\Box$

\section{Good splittable algebras}

In this section we prove a decomposition theorem for good splittable division 
algebras. This theorem shows how the studying of good splittable division 
algebras can be reduced to the studying of division algebras with  simple 
described structure. So, good splittable algebras are the most easy and good 
algebras to study. 

\begin{lemma}
\label{goodspl}
Let $D$ be a good splittable division algebra, $F=Z(D)$, and let 
$Z(\bar{D})=\bar{F}(s)$ be a purely inseparable over $\bar{F}$ field of degree 
$p=char D>0$. Let $u:\bar{D}\hookrightarrow D$ be a good embedding.

Then there exists a parameter $z$ such that 
${}^{(z,u)}_{-i}\delta_j=0$ for $j>i$, where $i=i(z,u)$ is a local height, and 
$u({}^{(z,u)}\delta_i(s))=x$, 
where $x\in Z(D)$. Moreover, $(i,p)=1$.
\end{lemma}

{\bf Proof.} Since $Z(\bar{D})/\bar{F}$ is  a purely inseparable extension, 
${}^{(z)}\alpha |_{Z(\bar{D})}=id$ for any parameter $z$.  By Skolem-Noether 
theorem there exists a parameter $z$ in $D$ such that ${}^{(z)}\alpha =id$. 
Suppose 
${}^{(z,u)}\delta_i(s)=0$, where $i=i(z,u)$. Then 
${}^{(z,u)}\delta_i|_{Z(\bar{D})}=0$, since $u$ is a good embedding and 
$Z(\bar{D})/\bar{F}$ is a simple extension. So, ${}^{(z,u)}\delta_i$ is an 
inner derivation by Scolem-Noether theorem, and by lemma \ref{ozamene}, (i) 
there exists a parameter $z'$ such that ${}^{(z',u)}\delta_i=0$, 
${}^{(z')}\alpha =id$. 

So, we can assume ${}^{(z,u)}\delta_i(s)\ne 0$ for some parameter $z$. Since 
$s^p\in Z(D)$, by lemma \ref{ppp} we have $(i,p)=1$. 
Since ${}^{(z,u)}\delta_i$ is a derivation, 
${}^{(z,u)}\delta_i(s)\in Z(\bar{D})$ (see the arguments in lemma \ref{(5)}, (ii)). 
Since $(i,p)=1$, there exists $k$ such that $p| (1-ki)$. So, 
by lemma \ref{ozamene}, (iii), for the parameter 
$z'=({}^{(z,u)}\delta_i(s))^k$ we have ${}^{(z')}\alpha =id$, 
${}^{(z',u)}\delta_i(s)\in \bar{F}$, i.e. $u({}^{(z',u)}\delta_i(s))\in Z(D)$. 
Since $s^p\in Z(D)$, by lemma \ref{vtorinv} we must have $d(u,s)=\infty$. In 
the proof of lemma \ref{vtorinv}, (i) was shown that $d(u,s)=d'(u,z,s)$ for 
some parameter $z$, and the construction of this element uses lemma 
\ref{ozamene}, (ii), so it preserves the initial values of ${}^{(z')}\alpha$, 
${}^{(z',u)}\delta_i$. So, ${}^{(z,u)}_{-i}\delta_j=0$ for $j>i$ and the lemma 
is proved.\\
$\Box$

\begin{prop} 
\label{razlozhenie}
Let $D$ be a splittable division algebra. Then
we have
$D\cong D_1\otimes_F D_2$, where $D_1, D_2$ are splittable division algebras 
such that
$D_1$ is an inertially split algebra. 

If $D$ is a good splittable division algebra, then 
 $Z(\bar{D_2})/\bar{F}$ is
a purely inseparable extension and $D_2$ is a good splittable algebra ($D_1$ 
or $D_2$ may be trivial).

So, $D\sim A\otimes_FB\otimes_FD_2$, where $A$ is a cyclic division algebra
and $B$ is an unramified division algebra.
\end{prop}

{\bf Proof.} If $char D=0$, the proposition is obvious, so we assume $char 
D>0$.

By \cite{P}, p.261, $D\cong D_1\otimes_F\ldots \otimes_FD_k$,
where $[D:F]\eq p_1^{r_1}\ldots p_k^{r_k}$ and $[D_i:F]\eq p_i^{r_i}$. Let
$p_2\eq p$. Since $D_i$ are defectless over $F$, $D_1,D_3,\ldots
D_k$ are inertially split. Therefore, by theorem \ref{Cohen} the algebra 
$B=D_1\otimes D_3\otimes
\ldots \otimes D_k$ is good splittable.

Assume first that $D$ is good splittable.  
By proposition 1.7. in \cite{JW}, if $s\in Z(\bar{D})$ is an element such that 
$\alpha (s)=s$, then this element is a purely inseparable element over 
$\bar{F}$. So, if $D$ is a good splittable division algebra, then 
 by lemma \ref{ppp} $D_2$ is either inertially split or $Z(\bar{D_2})/\bar{F}$ 
is a purely inseparable extension. For, otherwise there exists an element 
$s\in Z(\bar{D_2})\subset Z(\bar{D})$ as above and by proposition \ref{X} $p|i(u,s)$ for any  embedding $u$. If $u$ is a good embedding, then $s^{p^k}\in Z(D)$ for some $k$, a contradiction. 

So, we assume below $Z(\bar{D_2})/\bar{F}$ is a purely inseparable extension.
Now, we have (see, e.g. th.1 in \cite{Mor}) $\bar{D}\cong 
\bar{D_2}\otimes_{\bar{F}}\bar{B}$ and so $u(\bar{D})\cong 
u(\bar{D_2})\otimes_{u(\bar{F})}u(\bar{B})$, where $u$ is a good embedding. 
So, $E=u(Z(\bar{D_2}))$ is a purely inseparable field over $u(\bar{F})\subset 
Z(D)$.

Consider the field $E'=u(K)\otimes_{u(\bar{F})}F$, where $K$ is a maximal 
separable subfield in $\bar{B}$. This is an inertial lift of $K$ in $D$. 
Consider the centralizer 
$C=C_D(E')\cong D_2\otimes_FE'$. Let $M$ be a maximal subfield in $\bar{D_2}$. 
Note that $u(\bar{D_2})\subset C$, so $L\subset C$, where  
$L=u(M)F$ is the composit of  $u(M)$ and $F$, and $E\subset L$. 
Note that $[L:F]=ind D_2=ind C$. The field $L$ splits $C$ by dimension 
arguments. So, it must split $D_2$, since $([E':F],p)=1$, and $D_2$ is a 
$p$-algebra. Therefore, $L$ is isomorphic to a maximal subfield in $D_2$, so 
$D_2$ contain a copy of purely inseparable "unramified" subfield, whose 
residue field is isomorphic to $Z(\bar{D_2})$. Therefore, 
$D_2$ is a god splittable algebra. For, the centralizer of this field is an 
unramified division algebra, so by theorem \ref{Cohen} is splittable. So, 
$D_2$ is good splittable if the purely inseparable field is good splittable. 
But it is good splittable since it contains a subfield isomorphic to 
$u(Z(\bar{D_2}))$ by the construction. (Another way to see it is to use 
arguments from lemma \ref{simple} to show that there exists an appropriate 
$p$-basis).   

Let $D$ be a splittable algebra. Then the same arguments as in the previous 
paragraph show that $L$ is isomorphic to a maximal subfield in $D_2$ (it is 
not important that $Z(\bar{D_2})/\bar{F}$ may be not a purely inseparable 
extension).  
Now, the composit $EF\subset L$, $EF\ne L$, since every element from $E$ 
commute with $u(\bar{D_2})$, where $u$ is some fixed embedding.  So we must 
have $\overline{C_{D_2}(EF)=\bar{D_2}}$ and $C_{D_2}(EF)$ is an unramified 
division algebra. Therefore, $D_2$ is splittable division algebra.

Decomposition theorems \cite{JW}, Thm. 5.6-5.15
complete the proof.\\
$\Box$

This proposition shows that the study of splittable division algebras can be 
reduced to the study of splittable $p$-algebras. So, below in this section and 
in the next section we will deal with $p$-algebras only. 

\begin{prop}
\label{555}
Let $D$ be a good splittable division algebra such that 
$Z(\bar{D})/\overline{Z(D)}$ is
a purely inseparable extension. Then $D\cong D_1\otimes_{Z(D)}D_2$, where
$D_1$ is an unramified division algebra and $D_2$ is a good splittable 
division algebra such
that $\bar{D_2}$ is a field, $\bar{D_2}/\overline{Z(D)}$ is a purely
inseparable extension, $[\bar{D_2}:\overline{Z(D)}]\eq
[\Gamma_{D_2}:\Gamma_{Z(D)}]$.
\end{prop}

{\bf Proof.} The proof is by induction on the degree 
$[Z(\bar{D}):\overline{Z(D)}]$. 

Assume $[Z(\bar{D}):\overline{Z(D)}]=p$. Let ${}^{(z,u)}\delta_i$ be the map 
from lemma \ref{goodspl}. Then ${}^{(z,u)}\delta_i^p$ is a derivation trivial 
on the centre $Z(D)$, hence by Scolem-Noether theorem it is an inner 
derivation.

We claim that $z^p\in Z(D)$. 
We have 
$$
z^{-i}az^i\eq a+{}_{-i}\delta_i(a)z^i, \mbox{\quad} a\in u(\bar{D})
$$ 
Therefore,
$$
z^{-pi}az^{pi}\eq a+{}_{-i}\delta_i^p(a)z^{pi}, \mbox{\quad} a\in u(\bar{D})
$$
and 
$$
z^{pi}az^{-pi}\eq a+\delta'_1(a)z^{pi}+{\delta'_1}^2(a)z^{2pi}+\ldots ,
$$
where $\delta'_1\eq (-1){}_{-i}\delta_i^p=i^p\delta_i^p$. So,
$$
z^paz^{-p}\eq
a+\frac{1}{i}\delta'_1(a)z^{pi}+c_2\frac{1}{i^2}{\delta'_1}^2(a)z^{2pi}+\ldots
, 
$$
where $c_k$ are given by (\ref{(188)}) in lemma \ref{svva}. So, $z^p\in
Z(D)$ iff $\delta_i^p\eq 0$. Suppose $\delta_i^p\ne 0$. Consider an element
$Y\in Z(D)$, $w(Y)>0$. Let 
$$
Y\eq a_1z^p+\ldots , \mbox{\quad } a_1\in u(\bar{D}).
$$
First note that 
$$
Y\eq a_1z^p+a_2z^{2p}+a_3z^{3p}+\ldots , \mbox{\quad} a_i\in u(\bar{D})
$$
Indeed, $Y$ must satisfy $[Y,s]\eq 0$, where $s$ is a generator of 
$u(Z(\bar{D}))$ over $u(\bar{F})$. Since $s\in u(Z(\bar{D}))$ and 
$w([z^k,s])=k+i$ if $(k,p)=1$ and $w([z^k,s])=\infty$ otherwise, we then have
$[z^{i_k},s]\eq 0$ for every $k$, where 
$$
Y\eq \sum_{k\eq 1}^{\infty}a_kz^{i_k}
$$  
Therefore, $p|i_k$. 

Then, $Y$ must satisfy $Ya\eq aY$ for any $a\in u(\bar{D})$. Therefore,
$a_1,\ldots a_i\in u(Z(\bar{D}))$ and we must have 
$$
aa_{i+1}-a_{i+1}a\eq a_1\delta'_1(a)/i
$$
and
$$
aa_{2i+1}-a_{2i+1}a\eq a_i\delta'_1(a)+a_1c_2{\delta'_1}^2(a).
$$
Since $\Delta (a)\eq aa_{2i+1}-a_{2i+1}a$ is an inner derivation, we get  
${\delta'_1}^2\eq \delta$, where $\delta$ is a derivation, which is a
contradiction if $\delta\ne 0$ and $char D\ne 2$. In the last case we can use 
the same arguments with 
$a_{3i+1}$.  Therefore, ${\delta'_1}^2\eq \delta\eq 0$ and $\delta'_1\eq 0$,
and $z^p\in Z(D)$. 

Consider the algebra $W=u(Z(\bar{D}))((z))$. Since $z^p\in Z(D)$ and 
$u(\bar{F})\subset Z(D)$, we have $Z(W)=u(\bar{F})((z^p))=F$. So, $D\cong 
W\otimes_FC_D(W)$ by Double Centralizer theorem. It is clear that $C_D(W)$ is 
an unramified division algebra. 

Now suppose the proposition is proved for 
$[Z(\bar{D}):\overline{Z(D)}]=p^{k-1}$. By Albert's theorem (th.13 in 
\cite{Al}) $D_2$ then is a cyclic algebra as a product of cyclic subalgebras 
$D_i$, where $\bar{D_i}/\bar{F}$ is a simple purely inseparable extension and 
$D_i$ is a good splittable algebra. 

Assume $[Z(\bar{D}):\overline{Z(D)}]=p^{k}$. 
For a good embedding there exists a lift $\tilde{K}$ of a subfield
$\overline{Z(D)}\subset K\subset Z(\bar{D})$ such that the extension
$K/\overline{Z(D)}$ has degree $p$, i.e. $\bar{\tilde{K}}\eq K$,
$\Gamma_{\tilde{K}}\eq \Gamma_{Z(D_2)}$, $u(K)\subset \tilde{K}$, 
$\tilde{K}/Z(D)$ is a purely inseparable extension of degree $p$. By the 
induction hypothesis the centralizer $C_D(\tilde{K})\cong 
A_1\otimes_{\tilde{K}}A_2$, where $A_2$ is a cyclic division algebra and 
$\bar{A_2}$ is a field. Note that $\bar{A_2}=Z(\bar{D})$. 

By theorem 6 in \cite{Al} we can assume $A_2=(L/\tilde{K},\sigma ,a)$, where 
$a$ generate $\tilde{K}$ over $Z(D)$. So, $A_2$ contains a maximal purely 
inseparable Kummer subfield $E=\tilde{K}(y)$ with $y^{p^{k-1}}=a$, so 
$E=Z(D)(y)$. By theorem 3 in \cite{Al} $L=L_0\times \tilde{K}$, where $L_0$ is 
cyclic of degree $p^{k-1}$ over $Z(D)$ and $yx_0=\sigma (x_0)y$, where $x_0\in 
L_0$.

Consider the centralizer $B=C_D(L_0)$. We claim $B\cong B_1\otimes_{L_0} B_2$, 
where $B_2$ is a cyclic division algebra of degree $p$ and $B_2$ contains 
$\tilde{K}$. 

Note that
$B$ contains $Z(D)(a)=\tilde{K}$ and $A_1$. If $\tilde{K}L_0=L$ is 
"unramified" over $L_0$, then we apply the arguments for the first step of our 
induction to the algebra $B$. By construction, $B_2$ then will contain $L$, so 
$\tilde{K}$. Suppose $L$ is totally ramified over $L_0$ and let $z$ be a 
parameter of $L$, i.e. an element with the least possible positive mean of 
valuation on $L$. Since $L$ is purely inseparable over $L_0$, $z^p$ is a 
parameter of $L_0$. 

We have $W:=C_B(L)= C_D(L)\cong A_1\otimes_{\tilde{K}}L$ is an unramified 
division algebra. Consider an embedding $u':\bar{L}=\bar{L_0}\hookrightarrow 
L_0$. As it was shown in the proof of theorem \ref{Cohen} there is a lift 
$\tilde{u'}$ of $u'$, $\tilde{u'}:
\bar{W}\hookrightarrow  W$. Now consider the subalgebra 
$W'=\tilde{u'}(\bar{W})((z^p))$. We have $Z(W')=u'(\bar{L})((z^p))=L_0$, so 
$W'$ is an unramified subalgebra in $B$. By Double Centralizer theorem, 
$B\cong W'\otimes_{L_0}C_B(W')$, where $C_B(W')$ is a division algebra of 
degree $p$ and contains $L_0(z)=L$, so it contains $\tilde{K}$ and it is 
cyclic by Albert's theorem (th.12 in \cite{Al}).   

Now we can word by word repeat the arguments in the proof of theorem 12 in 
\cite{Al} to show that there exists a cyclic Galois extension $L'$ of $L_0$ 
which is cyclic Galois over $Z(D)$, and $y$ acts as a Galois automorphism on 
$L'/Z(D)$ which generates $Gal(L'/Z(D))$. So, there is the cyclic subalgebra 
$D_2=(L'/Z(D), ad(y), y^{p^k})$ in $D$.  Note that $A_2\subset D_2$, and $A_2$ 
is known to be a good splittable algebra with 
$[\bar{A_2}:\overline{Z(A_2)}]=[\Gamma_{A_2}:\Gamma_{Z(A_2)}]$. Since 
$\bar{A_2}=\bar{D_2}$ and $Z(A_2)=\tilde{K}$ is a purely inseparable extension 
of $Z(D)$, $D_2$ is a good splittable algebra such that $\bar{D_2}$ a field 
and $[\bar{D_2}:\overline{Z(D)}]=[\Gamma_{D_2}:\Gamma_{Z(D)}]$. By Double 
Centralizer theorem $D\cong D_1\otimes_{Z(D)}D_2$, where $D_1=C_D(D_2)$ must 
be an unramified division algebra, which completes the proof.\\
$\Box$

Combining all results in this section, we get the following theorem.

\begin{th}
\label{itog}
Let $D$ be a finite dimensional good splittable central division algebra over 
a field $F=k((t))$.  

If $char ({F})\eq p>0$, 
then $D\cong D_1\otimes_{F}D_2\otimes_{F}A_1\otimes_{F}\ldots 
\otimes_{F}A_m$, 
where $A_i$ are cyclic division algebras such that 
$[\bar{A_i}:\overline{Z(D)}]\eq
[\Gamma_{A_i}:\Gamma_{Z(D)}]$ and $\bar{A_i}/\overline{Z(D)}$ are simple 
purely inseparable field extensions, $D_1$ is an inertially split division
algebra, $(ind (D_1), p)\eq 1$, $D_2$ is an unramified division algebra ($D_1, 
D_2, A_i$ may be trivial).

If $char F=0$, then $D$ is an inertially split division algebra.  

\end{th}

\section{Splittability and good splittability}

In this section we collect some assorted results about a relation between splittable and good splittable division algebras and about splittable division algebras. 
We consider here only division algebras with the following property: 
$Z(\bar{D})/\overline{Z(D)}$ is a simple extension.

\begin{prop}
\label{cyclisity}
Let $D$ be a  central division algebra over $F$ of $char D=p>0$ such that 
$Z(\bar{D})=\bar{D}$ 
and $[Z(\bar{D}):\bar{F}]=p$. 

Then $D$ is a splittable algebra and the local height $i=i(u,z)$ (in the 
situatuion when it is defined, i.e. when $\alpha =id$) does not depend on $u$ and $z$.
It is a good splittable algebra if $(i,p)=1$. If $p|i$, then there exists a parameter $z$ such that $z^p\in Z(D)$ and any "unramified" maximal subfield is 
cyclic Galois. 

So, in both cases $D$ is a cyclic division algebra of degree $p$.  
\end{prop}

{\bf Proof.} Since $\bar{D}/\bar{F}$ is a simple extension, we have 
$[\bar{D}:\bar{F}]=[\Gamma_D:\Gamma_F]$. Indeed, 
consider the fields $E=F(s)$ and $E'=F(z)$, where $s$ is any element such that 
$\bar{s}$ is a primitive element of the extension $\bar{D}/\bar{F}$ and $z$ is 
any parameter of $D$. Then $[\bar{D}:\bar{F}]\le [E:F]\le 
[D:F]^{1/2}=([\bar{D}:\bar{F}][\Gamma_D:\Gamma_F])^{1/2}$, so 
$[\bar{D}:\bar{F}]\le [\Gamma_D:\Gamma_F]$. From another hand side, 
$[\Gamma_D:\Gamma_F]\le [E':F]\le 
([\bar{D}:\bar{F}][\Gamma_D:\Gamma_F])^{1/2}$, so 
$[\bar{D}:\bar{F}]=[\Gamma_D:\Gamma_F]$. 
So, $D$ is splittable division algebra of degree $p$.

If $Z(\bar{D})/\bar{F}$ is a separable extension, then $D$ is a good 
splittable algebra by theorem \ref{Cohen}. So, we assume it is a purely 
inseparable extension, $Z(\bar{D})=\bar{F}(\bar{u})$. For any lift $u$ of the 
element $\bar{u}$ 
let $u$ be an embedding constructed in lemma 
\ref{simple}, i.e. ${}^{(z,u)}\delta_j$ is defined by the values 
${}^{(z,u)}\delta_j(u^k)$ for any $j$.  By corollary \ref{ozamene3} the local 
height $i(u,z)$ does not depend on $z$, and by lemma \ref{ozamene2} $i(u,z)$ does not 
depend on $u$. 
For arbitrary embedding $u'$,  
since ${}^{(z,u')}\delta_{i(u',z)}$ is a derivation and $\bar{D}/\bar{F}$ is a simple extension, ${}^{(z,u')}\delta_{i(u',z)}$ is completely defined by a value at $\bar{u}$.  
Therefore, $i(u',z)=w(zu'(\bar{u})z^{-1}-u'(\bar{u}))$ and $i(u',z)$ is completely defined by the lift $u'(\bar{u})$. But arbitrary lift of $\bar{u}$ defines an embedding, on which we have proved $i$ does not depend. So, $i(u,z)$ does not depend on $z$ and $u$.

Now assume $p|i$. 

Using lemma \ref{ozamene}, we can assume without loss of generality that 
${}^{(z,u)}\delta_j=0$ if $j$ is not divisible by $p$. 

Indeed, if ${}^{(z,u)}\delta_j\ne 0$, then we apply lemma \ref{ozamene}, (ii) 
to show that there exists a parameter $z_j$ such that 
${}^{(z_j,u)}\delta_j(u)=0$ and ${}^{(z_j,u)}\delta_k={}^{(z,u)}\delta_k$ for 
$k<j$, 
${}^{(z_j)}\alpha =id$. Since ${}^{(z_j,u)}\delta_j$ is a derivation by 
proposition 
\ref{flyii} and by induction (similar arguments was already used in the proof 
of proposition \ref{X}), 
and since it 
is defined by the values on $u^k$, so by the values on $u$, we have 
${}^{(z,u)}\delta_j= 0$. 
Since for $j_1>j_2$ we have $w(z_{j_1}-z_{j_2})> j_1-i$, the sequence 
$\{z_j\}$ convereges to a  parameter $z'$, which satisfies our condition. 

So, there exists the subalgebra $A=u(\bar{D})((z^p))$. Let's show that 
$Z(D)\subset A$. 
Note that every element $a\in D$ can be written as $a=a_0+a_1z+\ldots 
+a_{p-1}z^{p-1}$, where $a_i\in A$. Note that $z^kAz^{-k}\subset A$ for every 
$k$. So, if $a\in Z(D)$, then $za_jz^{-1}=a_j$ and $ua_jz^ju^{-1}=a_jz^j$ for 
every $j$. For $j>0$ we have $a_jz^j=\sum_k a_{jk}z^{kp+j}$, so by corollary 
\ref{ozamene3} $ua_jz^ju^{-1}\ne a_jz^j$. 
Therefore, $a=a_0\in A$. 

Since $A\ne D$, $A$ must be commutative, so $z^p\in Z(D)$. Moreover, $A/Z(D)$ 
is cyclic Galois. Since  the arguments work for arbitrary lift $u$ of the 
element $\bar{u}$, arbitrary "unramified" maximal subfield in $D$ must be 
Galois over $F$. 

Now let $(i,p)=1$. 

Using lemma \ref{ozamene}, (iii) we can find a parameter $z$ and a primitive 
element $s\in \bar{D}$ 
such that ${}^{(z,u)}\delta_i(s)=sc$, where $c\in \bar{F}$. Indeed, since 
$(i,p)=1$, there exists $k$ such that $1-ki$ is divisible by $p$. So, by lemma 
\ref{ozamene}, (iii) for a parameter 
$z'=u({}^{(z,u)}\delta_i(\bar{u})^k)z$ we have ${}^{(z',u)}\delta_i(\bar{u})\in 
\bar{F}$, so by lemma \ref{ozamene2}, (iii) 
${}^{(z',u)}\delta_i(s)=1$, where 
$s=\bar{u}{}^{(z',u)}\delta_i(\bar{u})^{-1}$. 
Now, there exists $k_1$ such that $-ik_1-1$ is divisible by $p$, so for 
$z''=s^{k_1}z'$ we have ${}^{(z'',u)}\delta_i(s)=sc$, where $c=s^{-ik_1-1}\in 
\bar{F}$. It is easy to see  that, since $s=\bar{u}a$, where $a\in 
\bar{F}$, the map  ${}^{(z,u)}\delta_j$ is uniquely defined also by 
${}^{(z,u)}\delta_j(s^k)$, so by ${}^{(z,u)}_m\delta_l(s)$ for $l\le j$. 
So, we assume without loss of generality that $s=\bar{u}$, $z=z''$. 

Using lemma \ref{ozamene2}, (ii) we can find a converge sequence $\{u_j\}$, 
$u_j\in D$, $j\ge i$ such that $u_{j+1}=u_j+b_jz^{j+1-i}$, $u_i=u$, $b_j\in 
u_j(\bar{D})$ (here $u_j$ is an embedding defined by $u_j$, see lemma 
\ref{simple}) and 
${}^{(z,u_j)}_m\delta_k(\bar{u})\bar{u}^{-1}\in \bar{F}$ for all $k\le j$ and 
all $m$. 

Indeed, suppose it is true for $j\ge i$. Let 
${}^{(z,u_j)}_m\delta_{j+1}(\bar{u})=a_0+\ldots a_{p-1}\bar{u}^{p-1}$, $a_k\in 
 \bar{F}$. Since 
${}^{(z,u_j)}_m\delta_i={}^{(z,u)}_m\delta_i=m{}^{(z,u)}\delta_i$, we have 
$$
{}^{(z,u_j)}_m\delta_i(a_k\bar{u}^k)-\frac{\partial }{\partial 
\bar{u}}({}^{(z,u_j)}_m\delta_i(\bar{u}))a_k\bar{u}^k= (k-1)mca_k\bar{u}^k.
$$
So, $u_{j+1}=u_j-u_j(\sum_{k, k\ne 
1}(k-1)^{-1}m^{-1}c^{-1}a_k\bar{u}^k)z^{j+1-i}$ will satisfy our condition. 

We will denote by $u$ now  a limit of the sequence $\{u_j\}$. Using induction 
and proposition \ref{flyii} one can easily show that 
${}^{(z,u)}_m\delta_j(\bar{u}^k)\bar{u}^{-k}\in \bar{F}$ for any integer $k$. 
So, there is the subalgebra $A=u(\bar{F})((z))$ in $D$. Using similar 
arguments as in the case $p|i$, one can show that $A$ contains $Z(D)$. Since 
$A\ne D$, it must be commutative, so $u^p\in Z(D)$. Then $u$ is a good 
embedding, which completes the proof. \\
$\Box$

Let $D$ be a splittable division algebra and let $Z(\bar{D})/\overline{Z(D)}$ 
be a purely inseparable extension. As it was shown in the proof of lemma 
\ref{goodspl}, then there exists a parameter $z$ in $D$ such that 
${}^{(z,u)}\delta_{i}|_{Z(\bar{D})}\ne 0$, where $i=i(u,z)$ is a local height. 
Though $D$  may be not a good splittable algebra, the arguments from there are 
valid for every splittable algebra.  We will call such a parameter {\it an 
appropriate parameter}, and the number $i(u)=\max_zi(u,z)=i(u,z)$  for an 
appropriate parameter {\it a semilocal height}. Let's prove the following 
simple lemma.

\begin{lemma}
\label{simple2}
Let $D$ be a splittable central division $p$-algebra over $F$, where  $p=char 
D>0$, and let $Z(\bar{D})=\bar{F}(s)$ be a simple extension over $\bar{F}$. 
Then 

i) there exists an embedding $u$ such that  
${}^{(z,u)}_l\delta_j|_{Z(\bar{D})}$ is defined by the values 
${}^{(z,u)}_l\delta_j(s^k)$ for any $j,l,z$ (as in lemma \ref{simple});

ii) $[Z(\bar{D}):\bar{F}]=[\Gamma_{D}:\Gamma_{F}]$;

iii) if $\alpha |_{Z(\bar{D})}\ne id$ or $i(u)$ is divisible by $p$, then 
there exists a subalgebra $A=u(\bar{D})((z))$ for some appropriate parameter 
$z$ such that $Z(D)\subset Z(A)$. Moreover, $Z(A)$ is a cyclic Galois 
extension over $Z(D)$.
\end{lemma}

{\bf Proof.} i) 
 For arbitrary embedding $u$ consider the field $E=u(Z(\bar{D}))F\subset D$ 
and the centralizer $W=C_D(E)$. We have $\bar{W}=\bar{D}$ and so 
$Z(\bar{W})=\bar{E}$. Therefore, $W$ must be an unramified division algebra, 
and 
by theorem \ref{Cohen} there exists a lift on $\bar{W}$ of arbitrary  
embedding 
$\bar{E}\hookrightarrow E$. Now we can take an embedding defined by the 
element $s$ as in 
lemma \ref{simple}. It's lift will be desired embedding. We will denote this 
embedding also by $s$. 

ii) By proposition 1.7. in \cite{JW} the basic homomorphism $\theta_D$ (see 
introduction) is surjective. So, it is sufficient to prove the assertion only 
for the centralizer $C_D(K)$, where $K$ is a lift of a Galois part of the 
extension  $Z(\bar{D})/\bar{F}$. So, we will assume below $Z(\bar{D})/\bar{F}$ 
is a purely inseparable extension. 

Consider a maximal separable subfield $M$ in $\bar{D}$, and let $M'$ be a 
separable part of the extension $M/\bar{F}$. By \cite{JW}, th.2.8, th.2.9. 
there exists an inertial lift of $M'$ in $D$, say $\tilde{M}$. Consider the 
centralizer $B=C_D(\tilde{M})$. Then $\bar{B}$ is a field. Our assertion will 
be proved if we show it for $B$, since $[\tilde{M}:F]=ind (\bar{D})$ and 
$[D:F]=ind(\bar{D})^2[Z(\bar{D}):\bar{F}][\Gamma_{D}:\Gamma_{F}]$. 

Since $\bar{B}/\overline{Z(B)}$ is a simple extension, we can repeat the 
arguments from the beginning of proposition \ref{cyclisity}. 

iii) If $\alpha |_{Z(\bar{D})}\ne id$, consider the parameter $z$ from 
proposition \ref{X}. Then, clearly, $A=u(\bar{D})((z))$ will be a subalgebra 
with the center $K$, which is an inertial lift of a Galois part of the 
extension  $Z(\bar{D})/\bar{F}$. 

Assume  $\alpha |_{Z(\bar{D})}= id$ and $i(u)$ is divisible by $p$. Let $z$ be 
an appropriate parameter. 
Using lemma \ref{ozamene}, we can prove that ${}^{(z,u)}\delta_j=0$ if $j$ is 
not divisible by $p$. 

Indeed, let ${}^{(z,u)}\delta_j\ne 0$ be the first map with this property for 
$(j,p)=1$. If ${}^{(z,u)}\delta_j|_{Z(\bar{D})}=0$, then we apply lemma 
\ref{ozamene}, (i) to 
show that there exists a parameter $z_j$ such that 
${}^{(z_j,u)}\delta_j=0$ and ${}^{(z_j,u)}\delta_k={}^{(z,u)}\delta_k$ for 
$k<j$, 
${}^{(z_j)}\alpha =id$, since  ${}^{(z,u)}\delta_j$ is a derivation by 
proposition 
\ref{flyii} and by induction (similar arguments was already used in the proof 
of proposition \ref{X}) and so it is an inner derivation by Scolem-Noether 
theorem. 

 If ${}^{(z,u)}\delta_j|_{Z(\bar{D})}\ne 0$, then we apply lemma 
\ref{ozamene}, (ii) to show that there exists a parameter $z_j$ such that 
${}^{(z_j,u)}\delta_j(s)=0$ and ${}^{(z_j,u)}\delta_k={}^{(z,u)}\delta_k$ for 
$k<j$, 
${}^{(z_j)}\alpha =id$. Since ${}^{(z_j,u)}\delta_j$ is a derivation 
and since 
its restriction on  $Z(\bar{D})$ 
is defined by the values on $s^k$, so by the values on $s$, we have 
${}^{(z,u)}\delta_j|_{Z(\bar{D})}= 0$, 
 and we reduce the problem to the previous case. 
Since for $j_1>j_2$ we have $w(z_{j_1}-z_{j_2})> j_1-i$, the sequence 
$\{z_j\}$ convereges to a  parameter $z'$, which satisfies our condition.

Therefore, there exists a subalgebra $A=u(\bar{D})((z'))$ in $D$. Using the 
same arguments as in proposition \ref{cyclisity} one can show that 
$Z(D)\subset Z(A)$ Since $z'$ preserves $A$, it preserves the centre $Z(A)$ 
>From the other hand side, it acts nontrivially on it. So, $Z(A)$ is a cyclic 
Galois extension of degree $p$, and $ad(z')$ generates its Galois group.  \\
$\Box$

This lemma shows that the study of splittable $p$-algebras over $F$ can be 
reduced to the study of splittable  $p$-algebras with a purely inseparable 
extension $Z(\bar{D})/\bar{F}$ and $(i(u),p)=1$. 

\begin{defi}
\label{ivariant}
Let $D$ be a splittable division $p$-algebra with a purely inseparable 
extension $Z(\bar{D})/\bar{F}$. 
For any element $a\in \bar{D}$ define the number
$$
d_D(a)= \max_{u,z} w(z^{-i(u,a)}u(a)z^{i(u,a)}- 
u(a)-u({}^{(z,u)}_{-i(u,a)}\delta_{i(u,a)}(a))z^{i(u,a)})\in \dn\cup\infty ,
$$
where parameters $z$ are taken from the set of appropriate parameters and 
$i(u,a)$ was defined in corollary \ref{ozamene3}.
\end{defi}

It seems that the number $d_D(a)$ will play the role of a higher order level 
in a splittable division algebra. We will see that it codes a part of 
information about a division algebra.

\begin{lemma}
\label{predvarit}
Let $D$ be a splittable division $p$-algebra, $p>2$, with a purely inseparable 
simple extension 
$Z(\bar{D})/\bar{F}$, let $u$ be some fixed embedding 
$u:\bar{D}\hookrightarrow D$. 

Suppose $Z(\bar{D})=\bar{F}(a)$ and 
$(i(u,a),p)=1$. Suppose $d(u,a)\le 2i(u,a)$. 

Let $z$ be a parameter such that ${}^{(z,u)}\delta_{i(u,a)}
({}^{(z,u)}_{-i(u,a)}\delta_{i(u,a)}(a))=0$, ${}^{(z)}\alpha =id$ and
${}^{(z,u)}_{-i(u,a)}\delta_{q}|_{\sdf_p(a)}=0$ for 
$i(u,a)<q<d(u,a)$. 
Put $j(k):=i(u,a^{p^k})$.

Suppose for every $k\ge 1$ 
a parameter $z_k$ such that 
${}^{(z_k,u)}_{-j(k)}\delta_{r}|_{\sdf_p(a^{p^k})}=0$ for 
$j(k)<r<d(u,a^{p^k})$ satisfy a condition   
${}^{(z_k,u)}\delta_{i(u,a)}={}^{(z,u)}\delta_{i(u,a)}$, ${}^{(z_k)}\alpha ={}^{(z)}\alpha$. 

Suppose for every 
$k\ge 1$ we have 
$d(u,a^{p^k})-j(k)=d(u,a)-j(0)$. 

Then 
the maps ${}^{(z,u)}_{w+(p-1-r)j(k)}\delta_{\zeta}$, $rj(k)<\zeta\le (r-1)j(k)+d(u,a^{p^k})$, 
$r\in \{1, \ldots , p-1\}$, $k\ge 0$  
satisfy the following properties: 
$$
{}^{(z,u)}_{w+(p-1-r)j(k)}\delta_{\zeta}|_{\sdf_p(a^{p^k})}=c_{w+(p-1-r)j(k),\zeta ,1}\delta +\ldots +
c_{w+(p-1-r)j(k),\zeta ,r}\delta^{r},
$$ 
where the 
derivation $\delta$ was defined in lemma \ref{(5)}, $c_{w+(p-1-r)j(k),\zeta ,r}\in Z(\bar{D})$,  
$c_{w+(p-1-r)j(k),\zeta ,r}\ne 0$ only if $\zeta = (r-1)j(k)+d(u,a^{p^k})$. 

Moreover,  $c_{w+(p-1-r)j(k), (r-1)j(k)+d(u,a^{p^k}),r}\ne 0$ if 
$w=i(u,a) \mbox{\quad mod\quad}p$; 
$$c_{w+(p-1-r)j(k),(r-1)j(k)+d(u,a^{p^k}), r}=r!c_{w+(p-r)j(k),(r-2)j(k)+d(u,a^{p^k}), r-1}{}_{w+(p-1-r)j(k)}\delta_{j(k)}(a^{p^k}),$$
and ${}^{(z_k,u)}_{w+(p-2)j(k)}\delta_{d(u,a^{p^k})}(a^{p^k})=
{}^{(z_k,u)}_{-j(k)}\delta_{d(u,a^{p^k})}(a^{p^k})$. 

\end{lemma}

{\bf Proof.} The proof is similar to the proof of lemma \ref{(5)}, (i). It is by induction on $r$ simultaneously for all $k\ge 0$.  

For $r=1$, using lemma \ref{triviall} and induction, one can easily show that 
${}^{(z_k,u)}\delta_{q}(a^{p^k})=-(j(k))^{-1}{}^{(z_k,u)}_{-j(k)}\delta_{q}(a^{p^k})$ for $j(k)\le q< d(u,a^{p^k})$ (we assume here $z_0=z$). 
By lemma \ref{vtorinv}, (i) we have $d(u,a)-i(u,a)=i(u,a) \mbox{\quad mod\quad}p$. So, by lemma \ref{vtorinv}, (ii) and by induction we have $j(k)=j(0) \mbox{\quad mod\quad}p$. 

So,  ${}^{(z_k,u)}_{w+(p-2)j(k)}\delta_{q}|_{\sdf_p(a^{p^k})}=0$ if $j(k)< q< d(u,a^{p^k})$ and 
${}^{(z_k,u)}_{w+(p-2)j(k)}\delta_{q}|_{\sdf_p(a^{p^k})}\ne 0$ only if $q=d(u,a^{p^k})$. 

Since ${}^{(z_k,u)}_{-j(k)}\delta_{j(k)}|_{\sdf_p(a^{p^k})}$ is a derivation and 
since, by proposition \ref{flyii}, (i), the map  
${}^{(z_k,u)}_{-j(k)}\delta_{d(u,a^{p^k})}|_{\sdf_p(a^{p^k})}$ must be a derivation, we have 
${}^{(z_k,u)}_{w+(p-2)j(k)}\delta_{d(u,a^{p^k})}(a^{p^k})\in Z(\bar{D})$. For, as it was shown in the proof of lemma \ref{(5)}, (ii) 
for any derivation $\delta$ we have $\delta (b)\in Z(\bar{D})$ for any $b\in Z(\bar{D})$. Since 
${}^{(z_k,u)}_{w+(p-2)j(k)}\delta_{d(u,a^{p^k})}(a^{p^k})= q_1{}^{(z_k,u)}_{-j(k)}\delta_{d(u,a^{p^k})}(a^{p^k})+
q_2{}^{(z_k,u)}_{m}\delta_{j(0)}({}^{(z_k,u)}_{-j(k)}\delta_{j(k)}(a^{p^k}))$ for some integer $q_1,q_2,m$, we have proved our assertion. 
So, $c_{w+(p-2)j(k),d(u,a^{p^k}), 1}\in Z(\bar{D})$.

If $w=j(0) \mbox{\quad mod\quad}p$, then 
 ${}^{(z_k,u)}_{w+(p-2)j(k)}\delta_{d(u,a^{p^k})}(a^{p^k})={}^{(z_k,u)}_{-j(k)}\delta_{d(u,a^{p^k})}(a^{p^k})$, since $w+(p-1)j(0)=0\mbox{\quad mod\quad}p$ and $char D> 2$. So, we have $c_{w+(p-2)j(k),d(u,a^{p^k}), 1}\ne 0$. 

Put now $t=a^{p^k}$. 
For arbitrary $r$ by proposition \ref{flyii}, (i) we have 
$$
{}^{(z_k,u)}_{w+(p-1-r)j(k)}\delta_{\zeta}(t^q)\eq
q{}_{w+(p-1-r)j(k)}\delta_{\zeta}(t)t^{q-1}+
$$
$$
{}_{w+(p-1-r)j(k)}\delta_{j(k)}(t)
\sum_{l\eq 0}^{q-2}{}_{w+(p-r)j(k)}\delta_{\zeta -j(k)}(t^{q-1-l})t^l+
$$
$$
{}_{w+(p-1-r)j(k)}\delta_{d(u,t)}(t)
\sum_{l\eq 0}^{q-2}{}_{w+(p-1-r)j(k)+d(u,t)}\delta_{\zeta -d(u,t)}(t^{q-1-l})t^l+
$$
$$
\sum_{i=d(u,t)+1}^{\zeta -1} {}_{w+(p-1-r)j(k)}\delta_{i}(t)
\sum_{l\eq 0}^{q-2}{}_{w+(p-1-r)j(k)+i}\delta_{\zeta -i}(t^{q-1-l})t^l.
$$
Using the same arguments as in the proof of lemma \ref{(5)},(i) we see that ${}_{w+(p-1-r)j(k)}\delta_{\zeta}(t^p)=0$ and 
${}_{w+(p-1-r)j(k)}\delta_{\zeta}|_{\sdf_p(t)}=c_{w+(p-1-r)j(k),\zeta ,1}\delta +\ldots +
c_{w+(p-1-r)j(k),\zeta ,p-1}\delta^{p-1}$. To show that $c_{w+(p-1-r)j(k),\zeta ,i}=0$ for $i>r$ it suffice, by formulae (\ref{recurrent}) in lemma \ref{(5)}, to show that all the maps in the formula above are represented in the form $c_1\delta +\ldots +c_{r-1}\delta^{r-1}$. Let us show it in details. 

Since $\zeta -d(u,t)-1<(r-1)j(k)$, by lemma \ref{(5)}, (ii) 
${}_{m}\delta_{\zeta -i}|_{\sdf_p(t)}=c_{m,\zeta -i,1}\delta +\ldots +c_{m,\zeta -i,r-2}\delta^{r-2}$ for any $i>d(u,t)$. 

If $w=j(0) \mbox{\quad mod\quad}p$, then 
$w+(p-1-r)j(k)+d(u,t)+(r-2)j(k)=0 \mbox{\quad mod\quad}p$. Since 
$\zeta -d(u,t)\le (r-1)j(k)$, by lemma \ref{(5)}, (ii) we have 
${}_{w+(p-1-r)j(k)+d(u,t)}\delta_{\zeta -d(u,t)}|_{\sdf_p(t)}=c_{w+(p-1-r)j(k)+d(u,t),\zeta -d(u,t),1}\delta +\ldots +c_{w+(p-1-r)j(k)+d(u,t),\zeta -d(u,t),r-2}\delta^{r-2}$. 

If $w\ne j(0) \mbox{\quad mod\quad}p$, then by the same reason we have 
${}_{w+(p-1-r)j(k)+d(u,t)}\delta_{\zeta -d(u,t)}|_{\sdf_p(t)}=c_{w+(p-1-r)j(k)+d(u,t),\zeta -d(u,t),1}\delta +\ldots +c_{w+(p-1-r)j(k)+d(u,t),\zeta -d(u,t),r-1}\delta^{r-1}$ and by lemma \ref{(5)}, (i) $c_{w+(p-1-r)j(k)+d(u,t),\zeta -d(u,t),r-1}\in Z(\bar{D})$ as a product of elements from $Z(\bar{D})$. 

At last, by the induction hypothesis  
${}_{w+(p-r)j(k)}\delta_{\zeta -j(k)}|_{\sdf_p(t)}=c_{w+(p-r)j(k),\zeta -j(k),1}\delta +\ldots +c_{w+(p-r)j(k),\zeta -j(k),r-1}\delta^{r-1}$ and $c_{w+(p-r)j(k),\zeta -j(k),r-1}\ne 0$ only if $\zeta -j(k)=(r-2)j(k)+d(u,t)$, and $c_{w+(p-r)j(k),\zeta -j(k),r-1}\in Z(\bar{D})$. Since ${}_{w+(p-1-r)j(k)}\delta_{j(k)}(t)\in Z(\bar{D})$, by formulae 
(\ref{recurrent}) we get $c_{w+(p-1-r)j(k),\zeta ,r}\in Z(\bar{D})$ and if $w=j(0) \mbox{\quad mod\quad}p$, then 
$c_{w+(p-1-r)j(k),\zeta ,r}\ne 0$ iff $\zeta =(r-1)j(k)+d(u,t)$, 
$$c_{w+(p-1-r)j(k), (r-1)j(k)+d(u,t),r}=r!c_{w+(p-r)j(k),(r-2)j(k)+d(u,t), r-1}{}_{w+(p-1-r)j(k)}\delta_{j(k)}(t)\ne 0.$$ 
The lemma is proved.\\
$\Box$

\begin{lemma}
\label{predvarit2}
Let $D$ be a division algebra as in lemma \ref{predvarit}. 
Suppose $d(u,a)\le 2i(u,a)$ and $char D>2$. 

Then for every $k$ there exists a parameter 
$z_k$ such that 
${}^{(z_k,u)}_{-j(k)}\delta_{r}|_{\sdf_p(a^{p^k})}=0$ for 
$j(k)<r<d(u,a^{p^k})$ and  ${}^{(z_k)}\alpha ={}^{(z)}\alpha$, 
${}^{(z_k,u)}\delta_{j(l)}={}^{(z,u)}\delta_{j(l)}$ for all $l\le k$ (we use here the notation defined in lemma  \ref{predvarit}).

Moreover, for every 
$k\ge 1$ we have 
$d(u,a^{p^k})-j(k)=d(u,a)-j(0)$ and 
$${}^{(z_{k},u)}_{-j(k)}\delta_{d(u,a^{p^k})}(a^{p^k})=-{}^{(z_{k-1},u)}_{-j(k-1)}
\delta_{d(u,a^{p^{k-1}})}(a^{p^k})c_{d(u,t)-j(k-1), j(k)-j(k-1),p-1},$$
where 
$c_{d(u,t)-j(k-1), j(k)-j(k-1),p-1}$ is defined in lemma \ref{predvarit}. 

\end{lemma}

{\bf Proof.} The proof is by induction on $k$. By lemma \ref{vtorinv} 
$d(u,a)=2j(0)\mbox{\quad mod \quad }p$ and $j(1)=d(u,a)+(p-1)j(0)$. So, by the induction hypothesis we can assume for arbitrary $k$ that $d(u,a^{p^{k-1}})=2j(0)\mbox{\quad mod \quad }p$ and $j(k-1)=j(0)\mbox{\quad mod \quad }p$, and $j(k)=d(u,a^{p^{k-1}})+(p-1)j(k-1)$. 

For the convinience we can start with a parameter $z=z_0$, which satisfy the conditions of lemma \ref{predvarit}. Indeed, taking an appropriate parameter $z$ and changing it by a parameter $u(c)z$ for an appropriate $c\in Z(\bar{D})$ (as in the proof of proposition \ref{cyclisity}), we can assume that ${}^{(z,u)}_{-j(0)}\delta_{j(0)}(a)\in Z(\bar{D})^p$. Now, using arguments from the proof of lemma \ref{vtorinv}, (i), we can find such a parameter $z_0$. 

The idea of the proof is the following. We prove first that 
${}^{(z_{k-1},u)}_{-j(k-1)}\delta_{j(k)+d(u,a)-j(0)}(a^{p^k})\ne 0$. Then we prove that there exists a parameter $z_k$ such that  ${}^{(z_{k},u)}_{-j(k)}\delta_{\zeta}(a^{p^k})= 0$ for $j(k)<\zeta <j(k)+d(u,a)-j(0)$ and ${}^{(z_{k},u)}_{-j(k)}\delta_{j(k)+d(u,a)-j(0)}
(a^{p^k})\ne 0$. It will be shown that $z_k$ satisfy the conditions of lemma. 

So, assume $j(k)\le \zeta \le j(k)+d(u,a)-j(0)=j(k)+d(u,a^{p^{k-1}})-j(k-1)$. Put 
$t=a^{p^{k-1}}$. 
By proposition \ref{flyii}, (i) we have 
$$
{}^{(z_{k-1},u)}_{-j(k-1)}\delta_{\zeta}(t^p)\eq
$$
$$
{}^{(z_{k-1},u)}_{-j(k-1)}\delta_{d(u,t)}(t)
\sum_{l\eq 0}^{p-2}{}^{(z_{k-1},u)}_{d(u,t)-j(k-1)}\delta_{\zeta -d(u,t)}(t^{p-1-l})t^l+\ldots +
$$
$$
{}^{(z_{k-1},u)}_{-j(k-1)}\delta_{\zeta -(p-1)j(k-1)}(t)
\sum_{l\eq 0}^{p-2}{}^{(z_{k-1},u)}_{\zeta -pj(k-1)}\delta_{(p-1)j(k-1)}(t^{p-1-l})t^l+
$$
$$
\sum_{i=\zeta -(p-1)j(k-1)+1}^{\zeta -1} {}^{(z_{k-1},u)}_{-j(k-1)}\delta_{i}(t)
\sum_{l\eq 0}^{p-2}{}^{(z_{k-1},u)}_{i-j(k-1)}\delta_{\zeta -i}(t^{q-1-l})t^l.
$$
By lemma \ref{(5)}, (i)  in the last sum 
${}^{(z_{k-1},u)}_{i-j(k-1)}\delta_{\zeta -i}|_{\sdf_p(t)}=c_{i-j(k-1),\zeta -i,1}\delta +\ldots +c_{i-j(k-1),\zeta -i,p-2}\delta^{p-2}$, since $\zeta -i<(p-1)j(k-1)$. So, this sum is equal to zero. 

By lemma \ref{(5)}, (ii) we have ${}^{(z_{k-1},u)}_{\zeta -pj(k-1)}
\delta_{(p-1)j(k-1)}|_{\sdf_p(t)}= c_{\zeta -pj(k-1),(p-1)j(k-1),1}\delta +\ldots +
c_{\zeta -pj(k-1),(p-1)j(k-1),p-1}\delta^{p-1}$ and 
$c_{\zeta -pj(k-1),(p-1)j(k-1),p-1}\ne 0$ iff $\zeta =j(k-1)=j(0)\mbox{\quad mod\quad}p$. 

By lemma \ref{(5)}, (i) we have 
${}^{(z_{k-1},u)}_{m}\delta_{q}|_{\sdf_p(t)}=c_{m,q,1}\delta +\ldots +c_{m,q,p-1}\delta^{p-1}$ for $(p-1)j(k-1)<q<(p-1)j(k-1)+d(u,a)-j(0)$, and by lemma \ref{predvarit} $c_{m,q,p-1}=0$. 
By lemma \ref{predvarit} we have 
${}^{(z_{k-1},u)}_{d(u,t)-j(k-1)}\delta_{\zeta -d(u,t)}|_{\sdf_p(t)}=
c_{d(u,t)-j(k-1), \zeta -d(u,t), 1}\delta +\ldots + c_{d(u,t)-j(k-1), \zeta -d(u,t), p-1}\delta^{p-1}$ with $c_{d(u,t)-j(k-1), \zeta -d(u,t), p-1}\ne 0$ if $\zeta -d(u,t)=j(0)$. 

So, we have the following picture: ${}^{(z_{k-1},u)}_{-j(k-1)}\delta_{\zeta}(t^p)\ne 0$ only if $\zeta =j(0)\mbox{\quad mod \quad }p$ or if $\zeta =j(k)+d(u,a)-j(0)$. In the last case $${}^{(z_{k-1},u)}_{-j(k-1)}\delta_{\zeta}(t^p)=-{}^{(z_{k-1},u)}_{-j(k-1)}
\delta_{d(u,t)}(t^p)c_{d(u,t)-j(k-1), j(k)-j(k-1),p-1},$$
 where 
$c_{d(u,t)-j(k-1), j(k)-j(k-1),p-1}$ can be calculated using lemma \ref{predvarit}. 

Let's show that there exists a parameter $z_{k}$ such that 
${}^{(z_{k},u)}_{-j(k-1)}\delta_{\zeta}(t^p)=0$ for $j(k)<\zeta <j(k)+ d(u,a)-j(0)$. 
By lemma \ref{ozamene}, (ii) there exists a change of parameters $z_{k-1}\mapsto z'=z_{k-1}+bz_{k-1}^{p+1}$ such that ${}^{(z',u)}_{-j(k-1)}\delta_{j(k)+p}(t^p)=0$. 
It suffice to prove that any such a change of  parameters as in lemma \ref{ozamene}, (ii) with $p|q$ changes only the values of maps ${}_{-j(k-1)}\delta_{\zeta}$ with 
$\zeta =j(0)\mbox{\quad mod \quad }p$. For, if it is true, we can make several changes and kill all nonzero maps ${}^{(z_{k-1},u)}_{-j(k-1)}\delta_{\zeta}$ with $j(k)<\zeta <j(k)+ d(u,a)-j(0)$, since they are derivations and therefore are completely defined by their values at $t^p$. 

To prove it, we can use the calculations in the proof of lemma \ref{ozamene}, (ii). Since 
$d(u,a)-j(0)\le j(0)$, it is easy to see that for a change $z\mapsto z'=z+bz^{kp+1}$, $p>2$ we have there 
$$
z'^{-j(k-1)}t^pz'^{j(k-1)}=t^p+{}^{(z,u)}_{-j(k-1)}\delta_{j(k)}(t^p)z^{j(k)}+\ldots +
{}^{(z,u)}_{-j(k-1)}\delta_{j(k)+j(0)}(t^p)z^{j(k)+j(0)}+\ldots .
$$
Since $z'=z+bz^{kp+1}$, any power $z^{l}$ can be expressed as a series in $z'$, all powers of which are equal to $l$ modulo $p$. So, this change will change only maps with right indexes equal to $j(k)$ modulo $p$. 
Since ${}^{(z_{k-1},u)}_{-j(k-1)}\delta_{\zeta}(t^p)\ne 0$ only if $\zeta =j(0)\mbox{\quad mod \quad }p$ for $\zeta <j(k)+d(u,a)-j(0)$, our assertion is proved. 

So, there exists a parameter $z_k$ we have: ${}^{(z_k,u)}_{-j(k-1)}\delta_{\zeta}(t^p)\ne 0$ only if $\zeta =j(k)+d(u,a)-j(0)$ or $\zeta =j(k)$. Since $z_k$ was constructed as a sequence of changes as in lemma \ref{ozamene}, (ii), we have ${}^{(z_k)}\alpha ={}^{(z_{k-1})}\alpha$ and 
${}^{(z_k,u)}\delta_{j(q)}={}^{(z_{k-1},u)}\delta_{j(q)}$ for any $q\le k$. 

At last, let's prove that ${}^{(z_k,u)}_{-j(k)}\delta_{\zeta}(t^p)\ne 0$ only if $\zeta =j(k)+d(u,a)-j(0)$ or $\zeta =j(k)$. But this follows immediately from the definition of these maps, since $j(k)=j(k-1)\mbox{\quad mod \quad}p$, $d(u,a)-j(0)\le j(0)$ and $char D>2$. In particular, 
${}^{(z_k,u)}_{-j(k)}\delta_{j(k)}(t^p)={}^{(z_k,u)}_{-j(k-1)}\delta_{j(k)}(t^p)$, 
${}^{(z_k,u)}_{-j(k)}\delta_{j(k)+d(u,a)-j(0)}(t^p)={}^{(z_k,u)}_{-j(k-1)}\delta_{j(k)+d(u,a)-j(0)}(t^p)$. 

The lemma is proved.\\
$\Box$

Now we can prove the following theorem. 

\begin{th}
\label{posledn}
Let $D$ be a division $p$-algebra of $char D=p>2$ with the center $Z(D)=F$. Suppose $Z(\bar{D})=\bar{D}$ and $\bar{D}/\bar{F}$ is a simple purely inseparable extension, 
$\bar{D}=\bar{F}(a)$. Suppose that the semilocal height $i(u)$, which does not depend on the embedding $u$ in this case, is not divisible by $p$.  

Then $d_D(a)>i(u)$. 
\end{th}

{\bf Proof.} By lemma \ref{simple2}, (ii) $[\bar{D}:\bar{F}]=[\Gamma_D:\Gamma_F]$. So, the field $F(\tilde{a})$, where $\tilde{a}$ is a lift of $a$, is a maximal "unramified" subfield and therefore $D$ is a splittable division algebra. Obviously, $\alpha =id$. 

Since ${}^{(z,u)}\delta_{i(u,z)}$ is a derivation and $\bar{D}/\bar{F}$ is a simple extension, ${}^{(z,u)}\delta_{i(u,z)}$ is completely defined by a value at $a$. So, by lemma \ref{ozamene} $i(u,z)$ does not depend on $z$ and $i(u,z)=i(u)$. 
Therefore, $i(u)=w(zu(a)z^{-1}-u(a))$ and $i(u)$ is completely defined by the lift $u(a)$. From the other hand side, 
any lift $\tilde{a}$ of $a$ defines, by lemma \ref{simple}, an embedding $\tilde{a}$, and by lemma  
\ref{ozamene2} $i(\tilde{a})$ does not depend on $\tilde{a}$. So, $i(u)$ does not depend on $u$.  

The idea of the proof is following. We consider linear spaces which are the images of the maps ${}^{(z,u)}\delta_{j(k)}|_{\bar{F}(a^{p^k})}$ in $\bar{D}$, where $j(k)$ were defined in lemma \ref{predvarit2} and $z,u$ are fixed. We show that every such spase has zero intersection with each other if  $d_D(a)\le i(u)$. Then we show that this contradicts with the fact that $u(a)$ generate a finite dimensional space over $F$. 

So, assume $d_D(a)\le i(u)$. 
To calculate the spaces ${}^{(z,u)}\delta_{j(k)}({\bar{F}(a^{p^k})})\in \bar{D}$ we use lemmas \ref{vtorinv}, \ref{predvarit} and \ref{predvarit2}. We fix a parameter $z$ defined in lemma \ref{predvarit}. By lemmas \ref{simple}, \ref{ozamene2}, (iii) we can find a primitive element $\bar{u}\in \bar{D}$ of the extension $\bar{D}/\bar{F}$ such that ${}^{(z,u)}\delta_{j(0)}(\bar{u})=1$, where $u$ is an embedding defined in lemma \ref{simple} for some lift $u$ of the element $\bar{u}$. Using lemma \ref{ozamene}, (ii) we can find an embedding $u$ such that ${}^{(z,u)}\delta_{d(u,\bar{u})}(\bar{u})\notin {}^{(z,u)}\delta_{j(0)}(\bar{D})$. We fix this embedding. From lemmas \ref{ozamene}, \ref{ozamene2} immediately follows that $d(u,\bar{u})=d_D(\bar{u})=d_D(a)$. So, we assume without loss of generality $a=\bar{u}$. 

Put $J(k):={}^{(z,u)}\delta_{j(k)}(a^{p^k})$. Put $A(k):={}^{(z,u)}\delta_{j(k)}(\bar{F}(a^{p^k}))$, $A'(k):=\bar{F}(a^{p^{k+1}})\cdot a^{p^k(p-1)}J(k)$. 

We have $A(k)=\oplus_{q=0}^{p-2}\bar{F}(a^{p^{k+1}})\cdot a^{p^kq}J(k)$ and 
$\bar{D}\cdot J(k)=A(k)\oplus A'(k)$ as $\df_p$-linear spaces.  

From lemma \ref{vtorinv} follows that 
$${}^{(z,u)}\delta_{j(k)}(a^{p^k})={}^{(z_k,u)}\delta_{j(k)}(a^{p^k})=
q{}^{(z_{k-1},u)}_{-j(k-1)}\delta_{d(u,a^{p^{k-1}})}(a^{p^{k-1}})c_{d(u,a^{p^{k-1}})-j(k-1), (p-1)j(k-1), p-1},$$
where $q\in \df_p^*$, $z_k$ were defined in lemma \ref{predvarit}, $c_{d(u,a^{p^{k-1}})-j(k-1), (p-1)j(k-1), p-1}$ is calculated in lemma \ref{(5)}, (i) and it is not equal to zero by lemma \ref{(5)}, (ii), and ${}^{(z_{k-1},u)}_{-j(k-1)}\delta_{d(u,a^{p^{k-1}})}(a^{p^{k-1}})$ is calculated in lemma \ref{predvarit2}. 
By lemma \ref{predvarit2} we have ${}^{(z_{k-1},u)}_{-j(k-1)}\delta_{d(u,a^{p^{k-1}})}
(a^{p^{k-1}})=-j(k-1){}^{(z,u)}\delta_{d(u,a^{p^{k-1}})}(a^{p^{k-1}})$. 
Combining all these calculation together and using induction, we get  $J(k)=q_kJ(k-1)^pJ(1)=\tilde{q_k}J(1)^{p^{k-1}+p^{k-2}+\ldots +1}$ for $k\ge 1$, where $q_k\in \df_p$. 

Therefore, there is the following filtration 
$$
\bar{F}\subset \ldots \subset \bar{F}(a^{p^{k+1}})J(k+1)\subset \bar{F}(a^{p^k})J(k)\subset \ldots \subset \bar{D},
$$
and for every $k\ge 1$ we have $\bar{F}(a^{p^k})\cdot J(k)\subset A'(k-1)$. So, $A(k)\cap A(k_1)=\{0\}$ if $k\ne k_1$. 

Now consider an element $b\in F$ such that $\bar{b}=a^{p^l}$ for some $l>0$. We assume $l$ is a minimal possible integer. It exists, because $D$ is a finite dimensional algebra over $F$. Let 
$b=u(a^{p^l})+b_1z+\ldots $, where $b_k\in u(\bar{D})$. Put $I:=\min \{w(zb_kz^{k-1}-b_kz^k)\}$ (we assume here that $b_0=u(a^{p^l})$). Note that $I<\infty$, since  
by lemma \ref{predvarit2} $j(l)<\infty$, i.e. ${}^{(z,u)}\delta_{j(l)}(a^{p^l})\ne 0$. Now we must have 
$$
zbz^{-1}=\sum_{k=0}^{\infty}zb_kz^{k-1}= b+\sum_r {}^{(z,u)}\delta_{j(r)}(b_{q_r})z^I+\ldots =b, 
$$
where $b_{q_r}\in \bar{F}(a^{p^r})$ and $b_{q_r}\notin \bar{F}(a^{p^{r+1}})$. So, $\sum_r {}^{(z,u)}\delta_{j(r)}(b_{q_r})=0$, but it is impossible, since $A(k)\cap A(k_1)=\{0\}$ if $k\ne k_1$, a contradiction.

The theorem is proved. \\
$\Box$

{\bf Remark.} It would be interesting to know the answer on the following questions. 

i) Suppose $D$ is a division algebra as in the theorem \ref{posledn}. Does there exist a pair $(z,u)$ such that all nonzero maps ${}^{(z,u)}\delta_q$ satisfy the property $i(u)|q$? If it is true, there is a subalgebra $D'\subset D$ with $[D:D']<\infty$ and $D'$ has level 1 (see remark before lemma \ref{vtorinv}). So, we can reduce studying of $D$ to the algebra of level 1. 

ii) Is it true that $D$ is a good splittable algebra, i.e. cyclic? Probably, it is possible to apply our technique to give an answer to this question at least in the case of level 1.


\begin{thebibliography}{99}
\bibitem{Al}
{\it A.A. Albert} Normal division algebras of degree $p^l$ over $F$ of
characteristic $p$, Trans. Amer. Soc., 1936, 40, 112-126
\bibitem{AJ}
{\it R.Aravire, B.Jacob} $p$-algebras over maximally complete fields, Proc. Symp. Pure Math., vol.58, Part 2, AMS, Providence, 1995, pp.27-49
\bibitem{Az}
{\it G.Azumaya} On maximally central algebras, Nagoya J.Math. 2(1951), 119-150
\bibitem{Co}
{\it I. Cohen} On the structure and ideal theory of complete local rings,
Trans. Amer. Math. Soc., 59, 1946, p.54-106
\bibitem{C}
{\it P. M. Cohn} Algebra, vol.II, John Wiley \& Sons, 1991
\bibitem{JW}
{\it B.Jacob, A.Wadsworth} Division algebras over Henselian fields, J.Algebra,
128(1990), 126-179
\bibitem{DJ}
{\it A.J.de Jong} The period-index problem for the Brauer group of an algebraic surface, preprint. http://www-math.mit.edu/~dejong
\bibitem{Mor}
{\it P.Morandi} Henselisation of a valued division algebra, J.Algebra, 122,
1989, 232-243
\bibitem{P}
{\it R. S. Pierce} Associative Algebras, Springer-Verlag, New-York, 1982
\bibitem{PY}
{\it V.P.Platonov and V.I.Yanchevskii} Finite-dimensional skew fields, VINITI,
Moscow, 77(1991),
 144-262; English transl. in A.I.Kostrikin and
I.R.Shafarevich (eds.), Algebra IX, Encyclopaedia Math. Sci., vol. 77,
Springer, Berlin, 1995, pp. 121-239   
\bibitem{Sa}
{\it D. J. Saltman} Division algebras over discrete valued fields, Comm.
Algebra 8, 1980, 1749-1774
\bibitem{S} 
{\it O.F.S. Schilling} The theory of Valuations, Math.Surveys, Vol.4, Amer.
Math. Soc., Providence, RI, 1950
\bibitem{Ti}
{\it J.-P.Tignol} Algebres a division et extensions de corps sauvagement
ramifiees de degre premier, J.Reine Angew. Math. 404, 1990, 1-38
\bibitem{Zh}
{\it A.Zheglov} On a classification of two-dimensional skew fields,
Izvestiya RAN, 1, 2001,
pp. 30-65

 
\end{thebibliography}
\end{document}